\documentclass[12pt]{amsart}
\usepackage{amsfonts, latexsym, amssymb, amsgen, amsbsy, amstext, amsopn, amsmath}
\newcommand{\A}{\mathbb{A}}
\renewcommand{\H}{\mathbb{H}}
\newcommand{\C}{\mathbb{C}}

\newcommand{\G}{\mathbb{G}}

\newcommand{\M}{\mathbb{M}}
\newcommand{\N}{\mathbb{N}}

\newcommand{\R}{\mathbb{R}}

\newcommand{\V}{\mathbb{V}}

\newcommand{\cA}{\mathcal{A}}
\newcommand{\cB}{\mathcal{B}}

\newcommand{\cE}{\mathcal{E}}

\newcommand{\cG}{\mathcal{G}}

\newcommand{\cL}{\mathcal{L}}
\newcommand{\cN}{\mathcal{N}}

\newcommand{\cU}{\mathcal U}
\newcommand{\cM}{\mathcal{M}}
\newcommand{\cR}{\mathcal{R}}
\newcommand{\cT}{\mathcal{T}}

\newcommand{\cg}{\mathfrak{g}}
\newcommand{\ch}{\mathfrak{h}}
\newcommand{\ck}{\mathfrak{k}}
\newcommand{\cl}{\mathfrak{l}}

\newcommand{\cn}{\mathfrak{n}}

\newcommand{\cp}{\mathfrak{p}}
\newcommand{\cv}{\mathfrak{v}}
\newcommand{\cz}{\mathfrak{z}}

\newcommand{\Lip}{\mbox{\rm Lip}}

\renewcommand{\span}{\mbox{\rm span}}

\newcommand{\ep}{\varepsilon}

\newcommand{\ph}{\varphi}
\newcommand{\up}{\upsilon}
\newcommand{\sm}{\setminus}

\newcommand{\lls}{\mbox{\large $($}}
\newcommand{\rls}{\mbox{\large $)$}}
\newcommand{\lLs}{\mbox{\Large $($}}
\newcommand{\rLs}{\mbox{\Large $)$}}

\newcommand{\gope}{\circledcirc}

\newcommand{\ra}{\rightarrow}
\newcommand{\lra}{\longrightarrow}

\newcommand{\der}{\partial}
\renewcommand{\gope}{\circledcirc}

\newcommand{\cExp}{\mbox{$\cE\!$xp}}

\newtheorem{The}{Theorem}[section]
\newtheorem{Lem}[The]{Lemma}
\newtheorem{Def}[The]{Definition}
\newtheorem{Rem}[The]{Remark}
\newtheorem{Pro}[The]{Proposition}
\newtheorem{Cor}[The]{Corollary}
\newtheorem{Exa}[The]{Example}

\begin{document}

\title
%[Intrinsic graphs and $(\G,\R^k)$-regular sets in
%some H-type groups]
%{{\bf Intrinsic graphs and $(\G,\R^k)$-regular sets in
%some H-type groups}}
[Contact equations, Lipschitz extensions and isoperimetric inequalities]
{{\bf Contact equations, Lipschitz extensions and isoperimetric inequalities}}
\author{Valentino Magnani}
\address{Valentino Magnani, Dipartimento di Matematica \\
Largo Pontecorvo 5 \\ I-56127, Pisa}
\email{magnani@dm.unipi.it}

\maketitle
\begin{quote}
{\textsc{Abstract.}}
We characterize locally Lipschitz mappings and existence of
Lipschitz extensions through a first order nonlinear system of
PDEs. We extend this study to graded group-valued Lipschitz mappings
defined on compact Riemannian manifolds.
Through a simple application, we emphasize the connection between
these PDEs and the Rumin complex.
We introduce a class of 2-step groups, satisfying some abstract
geometric conditions and we show that Lipschitz mappings 
taking values in these groups and defined on subsets of the
plane admit Lipschitz extensions. We present several examples of
these groups, called Allcock groups, observing that their 
horizontal distribution may have any codimesion.
Finally, we show how these Lipschitz extensions theorems lead us to
quadratic isoperimetric inequalities in all Allcock groups.
\end{quote}

\tableofcontents

\newpage

\section{Introduction}

The local Lipschitz property for mappings between Euclidean spaces is characterized
by local $L^\infty$-bounds on distributional derivatives. If we replace the target
space even with the simplest sub-Riemannian manifold, as the Heisenberg group,
then the previous statement does not hold.
This elementary fact is the starting point of our study.

An important instance is the case of mappings defined on a Riemannian manifold,
that is related to the study of horizontal submanifolds.
In this connection, M. Gromov among other results
has treated various Lipschitz approximation theorems along with
Lipschitz extensions problems, see Section~3.5 of \cite{Gr1}.

It is well known that smooth Lipschitz mappings taking values into a
contact Carnot-Carath\'eodory manifold are horizontal and the converse
also holds, \cite{Gr1}. 
The main concern of this work is to understand this fact under Lipschitz regularity,
where both source and target spaces have the sub-Riemannian structure of
stratified groups. The horizontality condition for 
mappings of stratified groups yields the system of nonlinear first order PDEs \eqref{hcondF}, 
that only depends on the Lie algebra and that it can be written in terms of exponential coordinates. 

For general sub-Riemannian structures, these equations
can be replaced by the pull-back of forms defining the horizontal distribution,
according to Section~4.2 of \cite{Gr1}. Recall that in the case the target
of our mappings is Euclidean, then they are all horizontal and their
Lipschitz property is characterized in the usual way by bounds on
distributional horizontal derivatives. This fact holds for general 
Carnot-Carath\'eodory source spaces, \cite{GN1}.

Our point of view is that of considering mappings of stratified groups
as solutions to the system of equations \eqref{hcondF}.
Although the main results of this work are Theorem~\ref{CharLip} and Theorem~\ref{RMLip}, 
a substantial part of the paper is devoted to applications, regarding construction of 
both Lipschitz and non-Lipschitz mappings
and the relationship between Lipschitz extension theorems and isoperimetric inequalities.

The interest in geometric properties of mappings in the sub-Riemannian
setting has recently proved useful in connection with 
bi-Lipschitz embeddability of sub-Riemannian metric spaces
into different classes of infinite dimensional Banach spaces,
according to the remarkable work by J. Cheeger and B. Kleiner, \cite{CheKle}.

Our framework is that of {\em graded groups}, that are real, finite 
dimensional, connected, simply connected and nilpotent Lie groups,
with graded Lie algebra.
Notice that they might not be connected by rectifiable curves, according
to Example~\ref{GradParab}.
When the Lie subalgebra spanned by the first layer of the grading coincides 
with the whole algebra, we say that the group is {\em stratified}, \cite{FS}. 
The so-called {\em horizontal directions}
of the group are spanned by the left invariant vector fields belonging to the first
layer of the algebra, see Section~\ref{prel} for more details.
Stratified groups, also called Carnot groups,
represent the foremost models of nilpotent and simply connected
sub-Riemannian geometries. 
As we will see below, our techniques will allow for studying
mappings from a stratified group to a graded group.

The above mentioned horizontality for mappings of graded groups is well known as
{\em contact property}. This property corresponds to preserving horizontal directions. 
In Remark~\ref{cntpro}, we derive equations \eqref{hcondF} that correspond to this constraint.
We will refer to these equations as {\em contact equations}. 
For instance, every parametrization of a horizontal curve has the contact property 
by definition, hence the property of being horizontal can be equivalently stated in
terms of contact equations, that in this case become ODEs \eqref{eqhdiffcurve1}.
More generally, parametrizations of either Legendrian submanifolds in Heisenberg groups or
horizontal submanifolds in stratified groups are characterized by being
solutions to \eqref{hcondF}, where horizontal derivatives $X_j$
are replaced by usual partial derivatives.
Incidentally, these manifolds coincide with
$(\R^k,\M)$-regular sets, according to \cite{Mag11}.
Quasiconformal mappings of a stratified groups have the contact property, \cite{Pan2},
hence they solve the system \eqref{hcondF} a.e.
For instance, contact equations have been implicitly used in \cite{CapCow},
in relation with smoothness of 1-quasiconformal mappings between Carnot groups.
In connection with mappings with bounded distortion in
two step groups, contact equations explicitly appear in (3.2) of \cite{Dair}.
Certainly, many other interesting cases could be added from the existing literature.

In the present paper, we focus our attention on mappings defined on an open subset 
$\Omega$ of a stratified group $\G$, with graded group target $\M$. 
Lie algebras of $\G$ and $\M$ will be denoted by $\cG$ and $\cM$, respectively.
Turning to the initial question of characterizing the local Lipschitz property, 
we state the following
\begin{The}\label{CharLip}
Let $\Omega\subset\G$ be an open set and let 
$f:\Omega\lra\M$. Consider the mappings $F=\exp^{-1}\circ f$ and
$F_j=\pi_j\circ F$, hence $f=\exp\sum_{j=1}^\upsilon F_j$
and $F_j$ take values in $W_j$. If $F\in L^1_{loc}(\Omega,\cM)$,
then the following statements are equivalent:
\begin{itemize}
\item $f$ is locally Lipschitz,
\item
all the distributional derivatives $\nabla_{X_i}F$ belong to $L^\infty_{loc}(\Omega,\cM)$
and the system
\begin{equation}\label{hcondF}
\nabla_HF_j=\sum_{n=2}^\upsilon\frac{(-1)^n}{n!}\;\pi_j
\left([F,\nabla_HF]_{n-1}\right)
\end{equation}
$\mu$-a.e. holds in $\Omega$ for every $j=2,\ldots,\upsilon$,
\item
$f$ is a.e. P-differentiable and
$\nabla_HF=(\nabla_{X_1}F,\ldots,\nabla_{X_m}F)\in L^\infty_{loc}(\Omega,\cM)^m$,
\end{itemize}
where $(X_1,\ldots,X_m)$ is a basis of $V_1$ and $\mu$ is the Haar measure of $\G$.
The mappings $\pi_j:\cM\lra W_j$ indicate the canonical projection onto the 
$j$-th layer of the graded algebra $\cM$ of $\M$.
\end{The}
Under the stronger assumption that $\M$ is stratified, results of \cite{Pan2}
imply that the local Lipschitz property yields both the second and the third 
condition of Theorem~\ref{CharLip}.\-
Since a graded group need not be connected by rectifiable curves, then
the previously mentioned implications are obtained employing some technical tools
developed in \cite{Mag11}.

However, in this work we are interested in the converse to these implications, 
since we aim to achieve the Lipschitz property starting from the contact property.
Clearly, the standard smoothing argument used in Euclidean spaces to get
the Lipschitz property cannot be applied here. In fact, the mapping $F$ is a
weak solution to \eqref{hcondF}, but its mollification might no
longer be a solution of these equations, due to their nonlinearity. 
The direct use of Pansu differentiability does not seem to be of help, since this 
notion has no corresponding distributional version.
Here our point is that Theorem~\ref{CharLip} allows us to think of \eqref{hcondF} 
somehow as a distributional counterpart of Pansu differentiability.
Our approach boils down to making the problem a.e. one dimensional.
In fact, due to the distributional validity of \eqref{hcondF} and a Fubini's 
decomposition, one shows that for a.e. horizontal line in the domain,
the restriction of the mapping to this line is horizontal.
Then these restrictions are Lipschitz continuous with uniform Lipschitz constant
with respect to the homogeneous distance of the target.
This argument can be iterated for a basis of horizontal directions.
Since points are locally connected by piecewise horizontal lines
and any horizontal curve can be approximated by piecewise horizontal lines,
the local Lipschitz property follows.
Constructing globally Lipschitz functions on a domain $\Omega$ of a stratified group
clearly depends on the geometry of $\Omega$. In this respect,
every mapping that satisfies one of the equivalent conditions in Theorem~\ref{CharLip}
and also $\nabla_HF\in L^\infty(\Omega,\cM)^m$ is Lipschitz in $\Omega$,
if we assume that this set is a John domain.

The approach adopted to get Theorem~\ref{CharLip} also works replacing the open set 
of a stratified group with a complete Riemannian manifold, since one replaces 
horizontal lines with Riemannian geodesics, that are more manageable than 
sub-Riemannian's. By Theorem~\ref{LipContact} and
taking into account Remark~\ref{LipContactCl}, we are lead to the following
\begin{The}\label{RMLip}
Let $N$ be a compact connected Riemannian manifold and let $\nu$ be the canonical
Riemannian measure of $N$. Let $F:N\lra\cM$ be a Lipschitz mapping that 
$\nu$-a.e. satisfies
\begin{equation}\label{rkcont1}
\nabla F_j=\sum_{n=2}^\upsilon\frac{(-1)^n}{n!}\;\pi_j\left([F,\nabla F]_{n-1}\right)
\quad\mbox{for every $j=2,\ldots,\upsilon$.}
\end{equation}
Then $f:N\lra\M$, where $f=\exp\circ F$, is Lipschitz and there exists a geometric constant $C>0$ such that $\Lip(f)\leq C\;\Lip(F_1).$
\end{The}
A special instance of this theorem will be used in the proof
of Theorem~\ref{Allcock1Lipcon}, about the existence of Lipschitz extension.
Next, we discuss some applications of Theorem~\ref{CharLip}.
In fact, this theorem provides a general PDEs approach to construct either Lipschitz
or non-Lipschitz mappings. For instance, it is easy to construct smooth mappings 
in the Heisenberg group that are nowhere locally Lipschitz: it suffices 
to consider the parametrization of the vertical line in the
first Heisenberg group. In Subsection~\ref{hNnLip}, through contact equations,
we provide another case for a mapping of Heisenberg groups.
To construct examples of Lipschitz mappings, we use the natural relationship
between contact equations and the Rumin complex, \cite{Rum90}. 
In fact, an elementary computation only relying on contact equations and
using the complex property of Rumin differential allows us to determine all
smooth Lipschitz mappings of the Heisenberg group,
whose horizontal components are affine functions, see Subsection~\ref{RuminLinear}.
Our point here clearly is in the method, that rests on the solution of 
the simple contact equations specialized to this case.
It is also clear how our approach could be extended to
either higher dimensional Heisenberg groups or other variants.
Recall that affine Lipschitz mappings of the first Heisenberg group have been 
characterized in \cite{BHIT} through the explicit use of the Lipschitz condition.

Another application of Theorem~\ref{CharLip} concerns existence of Lipschitz 
extensions that can be interpreted as existence of solutions to \eqref{hcondF}
with assigned boundary datum. 
In general, this PDEs problem is very difficult, already for Euclidean spaces and
Heisenberg groups.
We will state and prove this characterization in Subsection~\ref{charextlip}.
This leads us to two possible methods to approach Lipschitz extension
problems in sub-Riemannian geometry: the geometric approach and the PDEs approach.

These two demanding projects cannot be treated here in depth, then we will
limit ourselves to show first elementary examples of how one can apply
both Theorem~\ref{CharLip} and Theorem~\ref{RMLip} to obtain Lipschitz extensions,
showing how this issue is essentially equivalent to finding isoperimetric inequalities.
We concentrate our attention on the geometric approach, following the method by D. Allcock to get quadratic isoperimetric inequalities in higher dimensional Heisenberg groups, \cite{All}.
Essentially, contact equations permit us to rephrase Allcock construction
in a larger class of two step groups that must satisfy some abstract geometric 
properties. We call these groups {\em Allcock groups}.
This leads us to a disk extension theorem in this class of groups, see 
Theorem~\ref{Allcock1Lipcon}.
Notice that M. Gromov pointed out how his disk extension theorem for contact simply 
connected compact Carnot-Carath\'eodory target gives a quadratic isoperimetric
inequality, see p.218 of \cite{Gr1}. As a final application of the tools developped
in this work,
we will also prove this implication for all Allcock groups, see Theorem~\ref{quadisop}.

In Section~\ref{allcockgroups}, we introduce Allcock groups and present several 
examples, that include higher dimensional quaternionic H-type groups,
the complexified Heisenberg group and other classes of two step groups.
It is easy to show that Allcock groups also include Heisenberg groups,
see Remark~\ref{AllcockHeisenberg}. 
Next, we state the following
%
%
%                         1-CONNECTEDNESS
%
%
\begin{The}[Disk extension theorem]\label{Allcock1Lipcon}
Let $\A\cl^n$ be an Allcock group, with $n\geq2$.
Then there exists a geometric constant $c>0$ such that
for every Lipschitz mapping $f:S^1\lra\A\cl^n$ there
exists a Lipschitz extension on the closed disk
$\tilde f:D\lra\A\cl^n$ such that
$\Lip(\tilde f)\leq c\;\Lip(f)$.
\end{The}
It is a rather general fact that Lipschitz extension theorems from spheres
to the corresponding higher dimensional disks
imply full Lipschitz extension theorems, from an Euclidean source space.
The argument of the proof essentially relies on Whitney cube
decomposition, according for instance to Theorem~1.2 of \cite{Almgren}.
A different statement of this fact can be found for mappings taking values 
in a compact Carnot-Carath\'eodory manifold, see p.219 of \cite{Gr1}. 
In the recent work by Lang and Schlichenmaier, \cite{LanSch},
using Nagata dimension, more general source spaces are considered. 
In Section~\ref{LipExtSect}, we use both terminology and results of this paper 
to show how Theorem~\ref{Allcock1Lipcon} yields the following
%
%
%                  LIPSCHITZ EXTENSIONS
%
%
\begin{Cor}\label{LipExtR2Aln}
Let $\A\cl^n$ be an Allcock group, with $n\geq2$.
Then the couple $(\R^2,\A\cl^n)$ has the Lipschitz extension property.
\end{Cor}
Section~\ref{Sectquadisop}
is devoted to quadratic isoperimetric inequalities for Allcock groups. 
We wish to stress that here we certainly do not exhibit novel techniques to get isoperimetric inequalities in nilpotent groups, but rather aim to clarify known
results, presenting some possibly new cases and eventually 
support future progress to study this vast and interesting issue 
by methods of sub-Riemannian Geometric Measure Theory.
\begin{The}[Quadratic isoperimetric inequality]\label{quadisop}
Let $\A\cl^n$ be an Allcock group, with $n\geq2$.
Then there exists a geometric constant $K>0$ such that 
for every Lipschitz loop $\Gamma:S^1\lra\A\cl^n$
there exists a Lipschitz map $f:D\lra\A\cl^n$ such that
$f_{|S^1}=\Gamma$ and the ``spanning disk'' $(D,f)$ satisfies the 
quadratic isoperimetric inequality
\begin{equation}\label{quadisoper}
{\mathcal H}_{\rho_0}^2\big(f(D)\big)\leq K\;\mbox{\rm length}_{\rho_0}(\Gamma)^2\,,
\end{equation}
where $\rho_0$ is a fixed Carnot-Carath\'eodory distance.
\end{The}
Recall that quadratic isoperimetric inequalities impose quadratic behaviour to the 
corresponding ``geometric Dehn functions". It is interesting to point out that
there is also a combinatorial notion of Dehn function, that is shown to be
equivalent to the geometric one, \cite{Bri}, \cite{BurTab}.
This represents a fascinating connection between Combinatorial Group Theory
and Geometric Group Theory. To have a glimpse of this vast research area,
we mention just a few references, \cite{Bri}, \cite{ECHLPT}, \cite{GrHypG}, 
\cite{GrAsympt}. 

In fact, quadratic isoperimetric inequalities in Heisenberg groups have also a
combinatorial proof, \cite{OlsSap}. More recently, R. Young have found various 
estimates for Dehn functions in some specials classes of stratified groups, \cite{You06}, \cite{You08}. 
In particular, for some central powers of two step groups he establishes quadratic isoperimetric inequalities using combinatorial methods, \cite{You06}.
It is interesting to compare his results with our Theorem~\ref{quadisop},
since the construction of central powers is somehow similar to the 
construction of an Allcock group $\A\cl_\cn^n$, starting from
its model algebra $\cn$, see Section~\ref{allcockgroups}.
 
Notice that the core of Allcock isoperimetric inequality is establishing 
this result in the ``symplectic part" $(\R^{2n},\omega)$ of the Heisenberg group
$\H^n$ using the Euclidean Hausdorff measure in $\R^{2n}$.
In this respect, he raises the question on the proper notion of area to read
isoperimetric inequalities in the Heisenberg group equipped with its Carnot-Carath\'eodory distance, see at p.230 of \cite{All}. 
Theorem~\ref{quadisop} answers this question showing that 
the 2-dimensional Hausdorff measure with respect to the Carnot-Carath\'eodory
distance works. This is a simple consequence of the sub-Riemannian
area formula, \cite{Mag}, that shows how the Euclidean surface measure of the 
projected surface in $\R^{2n}$ corresponds to its Hausdorff measure
with respect to the Carnot-Carath\'eodory distance in $\H^n$.
This works more generally for Allcock groups, see Proposition~\ref{srareapro}. 

In this case, we first obtain an Euclidean isoperimetric inequality in the ``multi-simplectic space'' $(\R^{mn},\omega)$, where the multi-symplectic form
$\omega$ is defined in \eqref{multsimpl}.
In fact, $\A\cl^n$ can be identified with $\R^{mn}\times\R^s$ with respect to 
suitable graded coordinates and the horizontal subspace given by $\R^{mn}$
inherits $\omega$ from the Lie algebra of $\A\cl^n$. The quadratic isoperimetric inequality in $(\R^{mn},\omega)$
is then the main point, corresponding to Theorem~\ref{projquadisoperthe}.
Then the sub-Riemannian area formula leads us to Theorem~\ref{quadisop}.

Since the Riemannian Heisenberg group is quasi-isometric to the sub-Riemannian one,
it is natural to expect the same isoperimetric inequality with respect to both 
distances, according to the case of finitely presented groups, \cite{Alo}. 
In fact, under the same assumtptions of Theorem~\ref{quadisop} we have 
\begin{equation}\label{riemquadisoper}
{\mathcal H}_\varrho^2\big(f(D)\big)\leq K\;\mbox{\rm length}_\varrho(\Gamma)^2\,,
\end{equation}
where $\varrho$ is the Riemannian distance obtained by the fixed left
invariant Riemannian metric defining the Carnot-Carath\'eodory
distance $\rho_0$, see Section~\ref{Sectquadisop}.
In fact, since $\rho_0$ is greater than or equal to $\varrho$ and
Proposition~\ref{lengthcurve} show that Riemannian
and sub-Riemannian lengths of a horizontal curve coincide,
it follows that \eqref{quadisoper} implies \eqref{riemquadisoper}.
Finally, we wish to remark that using the Euclidean Hausdorff measure
to study quadratic isoperimetric inequalities for instance in $\H^2$ does not seem promising, according to Example~\ref{exaheis}.
%
%
%
%\pagebreak

%%%%%%%%%%%%%%%%%%%%%%%%%%%%%%%%%%%%%%%%%%%%%%%%%%%%%%%%%%%%%%%%%%%%%%%%%%%%%%%%%
%
%
%
%
%               PRELIMINARIES AND KNOWN RESULTS
%
%
%
%
%%%%%%%%%%%%%%%%%%%%%%%%%%%%%%%%%%%%%%%%%%%%%%%%%%%%%%%%%%%%%%%%%%%%%%%%%%%%%%%%%

\section{Preliminaries and known results}\label{prel}

\subsection{Some elementary facts on graded groups}
A {\em graded group} is a real, finite dimensional,
connected and simply connected Lie group $\M$, whose Lie algebra $\cM$
can be written as the direct sum of subspaces $W_i$, called layers,
such that
\begin{eqnarray}\label{grading}
[W_i,W_j]\subset W_{i+j}
\end{eqnarray}
and $\cM=W_1\oplus\cdots\oplus W_\upsilon$.
The integer $\upsilon$ is the step of nilpotence of $\M$.
A graded group $\M$ is {\em stratified} if its layers satisfy
the stronger condition $[W_i,W_j]=W_{i+j}$, see for instance \cite{FS}.
The {\em horizontal tangent spaces} 
\[
H_x\M=\{Z(x)\mid Z\in W_1\}\subset T_x\M,\qquad x\in \M
\]
define all horizontal directions of the group, that are
collected into the so-called {\em horizontal subbundle} $H\M$.
We also define subbundles of higher order $H^j\M$, setting
\[
H^j_x\M=\{Z(x)\mid Z\in W_j\}\subset T_x\M,\qquad x\in \M
\]
We fix a norm $\|\cdot\|$ in $\cM$, then the bilinearity of
Lie brackets gives
\begin{equation}\label{lcst}
\|[X,Y]\|\leq\beta\, \|X\|\;\|Y\|\quad\mbox{for every}\quad X,Y\in\cM
\end{equation}
for some constant $\beta>0$, depending on the norm and on the algebra.
The grading of $\cM$ allows us to introduce a one-parameter group
of Lie algebra automorphisms $\delta_r:\cM\lra\cM$, defined as
$\delta_r(X)=r^i$ if $X\in V_i$, where $r>0$.
These mappings are called {\em dilations}.
Taking into account that the exponential mapping $\exp:\cM\lra\M$
is a diffeomorphism for simply connected nilpotent Lie groups,
we can read dilations in the group $\M$ through the mapping $\exp$
and mantain the same notation.
We fix a {\em homogeneous distance} $\rho$ on $\M$, namely, 
a left invariant continuous distance that is 1-homogeneous
with respect to dilations $\delta_r$.
We will use the convention $\rho(x)=\rho(x,e)$, where $e$ denotes
here the unit element of $\M$.
\begin{Exa}\label{GradParab}{\rm 
Let $\M=\R^n\times\R$ be equipped with the sum of vectors as commutative group operation
and define the parabolic distance $|(x,t)|=|x|+\sqrt{|t|}$ and dilations $\delta_r(x,t)=(rx,r^2t)$. Here we have the grading $V_1\oplus V_2$, where 
$V_1$ and $V_2$ can be identified with $\R^n\times\{0\}$ and $\{0\}\times\R$, respectively. 
Clearly, $\M$ is a 2-step graded group, but it is not stratified.}
\end{Exa}
Following the notation of \cite{Pan2}, Section~4.5, in the next definition
we introduce the iterated Lie bracket.
\begin{Def}{\rm
Let $X,Y\in\cM$. The $k$-th bracket is defined by
\begin{eqnarray}
[X,Y]_k=\underbrace{[X,[X,[\cdots,[X}_{\mbox{\tiny $k$ times}},Y],],\ldots,]
\quad\mbox{and}\quad[X,Y]_0=Y\,.
\end{eqnarray}
Notice that $[X,Y]_k=\mbox{ad}^k(Y)$.
}\end{Def}
We can express the group operation in the Lie algebra
by the Baker-Campbell-Hausdorff formula, that we present
in the following form
\begin{eqnarray}\label{absBCH}
X\gope Y=\sum_{j=1}^\upsilon	c_n(X,Y),
\end{eqnarray}
where $c_1(X,Y)=X+Y$, $c_2(X,Y)=[X,Y]/2$ and the subsequent
terms can be defined inductively, \cite{Vara}.
Notice that the sum \eqref{absBCH} has $\up$ addends,
since $X,Y$ belong to $\cM$ that has step of nilpotence $\upsilon$.
The addends $c_n$ in general are given by induction through
the following Baker-Campbell-Hausdorff-Dynkin formula
\begin{eqnarray}\label{kBCH}
&&(n+1)\,c_{n+1}(X,Y)=\frac{1}{2}\;[X-Y,c_n(X,Y)] \\
&&+\sum_{\substack{p\geq 1\\ 2p\leq n}}K_{2p}
\sum_{\substack{k_1,\ldots,k_{2p}>0 \\
k_1+\cdots k_{2p}=n}}[c_{k_1}(X,Y),[\cdots,[c_{k_{2p}}(X,Y),X+Y],],\ldots,],
\nonumber
\end{eqnarray}
see Lemma~2.15.3 of \cite{Vara}.
\begin{Lem}
Let $\nu>0$ and let $n=2,\ldots,\iota$. Then there exists
a constant $\alpha_n(\nu)$ only depending on $n$ and $\nu$ such that
\begin{eqnarray}\label{bilestim}
\|c_n(X,Y)\|\leq \alpha_n(\nu)\;\|[X,Y]\|
\end{eqnarray}
whenever $\|X\|,\|Y\|\leq\nu$.
\end{Lem}
{\sc Proof.}
Our statement is trivial for $n=2$, being $c_2(X,Y)=[X,Y]/2$.
Assume that it is true for every $j=2,\ldots,n$, with $n\geq2$.
We observe that $[c_{k_{2p}}(X,Y),X+Y]\neq0$ in \eqref{kBCH} implies $k_{2p}>1$, then inductive hypothesis yields
$$
\|c_{k_{2p}}(X,Y)\|\leq \alpha_{k_{2p}}(\nu)\; \|[X,Y]\|\,.
$$
Using this estimate in \eqref{kBCH} and observing that
$\|c_{k_i}(X,Y)\|\leq2\nu$, whenever $k_i=1$,
our claim follows. $\Box$
\begin{Def}{\rm
A {\em homogeneous subgroup} $H$ of $\G$ is a Lie subgrup that is closed
under dilations. Analogously, for subalgebras the same terminology is adopted.}
\end{Def}
\begin{Def}{\rm
Let $N$ and $H$ be homogeneous subgroups of $\G$, where
$N$ is normal, $N\cap H=\{e\}$ and $NH=\G$. Then $\G$ is an
inner semidirect product of $N$ and $H$ and we write
$\G=N\rtimes H$.}
\end{Def}
We denote by $H_X$ the one-dimensional subgroup of $\G$,
spanned by $\exp X$, where $X\in\cG$. 
Next, we recall a standard fact concerning direct sums of
homogeneous subalgebras, whose proof can be found for instance in \cite{Mag11}.
Recall that a homogeneous algebra $\cp$ satisfies $\delta_r\cp\subset\cp$
for every $r>0$.
\begin{Pro}\label{homdec}
Let $\cp$ and $\ch$ be homogeneous subalgebras of $\cG$
and let $P$ and $H$ denote their corresponding homogeneous subgroups,
respectively. Then the condition $\cp\oplus\ch=\cG$ is equivalent
to require that $P\cap H=\{e\}$ and $PH=\G$.
Furthermore, if one of these conditions hold, then the mapping
\begin{eqnarray}\label{phidecomp}
\phi:\cp\times\ch\lra\G,
\qquad \phi(W,Y)=\exp W\;\exp Y
\end{eqnarray}
is a diffeomophism.
\end{Pro}
\begin{Lem}\label{factorX}
Let $X\in V_1\sm\{0\}$. Then there exists a normal homogeneous
subgroup $N\subset\G$ such that $\G=N\rtimes H_X$.
\end{Lem}
{\sc Proof.}
Let $U_1\subset V_1$ be a subspace of dimension $\dim V_1-1$
that does not contain $X$ and set
$\cN=U_1\oplus V_2\oplus\cdots\oplus V_\iota$.
We notice that $\cN$ is a homogeneous ideal,
hence $N=\exp\cN$ is a normal homogeneous subgroup of $\G$.
By definition of $\cN$, we have $\cN\cap\span\{X\}=\{0\}$ and
this implies that $N\cap H_X=\{e\}$ and $NH_X=\G$, see for instance Proposition~\ref{homdec}. In other words $\G=N\rtimes H_X$. $\Box$

\subsection{Haar measure and Fubini's theorem.}
We consider a graded group $\G$ with grading $\cG=V_1\oplus\cdots\oplus V_\iota$
and we fix a left invariant Riemannian metric $g$ on $\G$.
Then the associated volume measure $\mbox{vol}_g$ is clearly left invariant
and defines the Haar measure of $\G$. We will denote by $\mu$ this
measure. Let $X\in V_1$ and $N$ be a homogeneous normal subgorup $N$ such that
$\G=N\rtimes H_X$. Then we have the following Fubini's theorem
with respect to this factorization.
\begin{Pro}\label{FubX}
Let $\mu$ be the Haar measure of $\G$.
Then for every measurable set $A\subset\G$, we have
\begin{eqnarray}\label{FubNX}
\mu(A)=\int_N\nu_X\left(A_n\right)\,d\mu_N(n)\,,
\end{eqnarray}
where $A_n=\{h\in H_X\mid nh\in A\}$. We have 
denoted by $\mu_N$ and $\nu_X$ the Haar measure of $N$ and of $H_X$, respectively,
\end{Pro}
{\sc Proof.}
We fix an orthonormal basis $(X_1,\ldots,X_q)$ of $\cG$
with respect to the metric $g$. In addition we assume that
this basis is adapted to the grading of $\cM$, such that $X_1$ is proportional to $X$.
By Proposition~\ref{homdec}, the mapping
\[
\psi:\R\times\R^{q-1}\lra\M, \qquad (\xi,t)\lra
\exp\Big(\sum_{j=2}^q x_jX_j\Big)\exp\big(x_1X_1\big)
\]
is a diffeomorphism. Our fixed basis also introduces the special
system of coordinates $F:\R^q\lra\G$, $F(x)=\exp\big(\sum_{i=1}^qx_iX_i)$
on $\G$. Now we observe that
\[
F_\sharp\cL^q=\mu,\quad \big(F_{|\{0\}\times\R^{q-1}}\big)_\sharp\cL^{q-1}=\mu_N,\quad
\mbox{ and }\quad \big(F_{\R\times|\{0\}}\big)_\sharp\cL^1=\nu_X\,,
\]
see for instance Proposition~2.3.47 of \cite{MagPhD}.
Thus, the fact that $X_1\in V_1$ implies that the mapping $F^{-1}\circ\psi$ has jacobian equal to one. Combining these facts with classical Fubini's theorem,
we get our claim. $\Box$

\subsection{Differentiability}

In the present subsection, we recall the notion of Pansu differentiability.
$\G$ and $\M$ denote two stratified groups and $\Omega$ is an open subset of $\G$.
\begin{Def}[h-homomorphism]{\rm
A group homomorphism $L:\G\lra \M$ such that
$L(\delta_r^\G x)=\delta_r^\M L(x)$ for every $x\in\G$ and $r>0$
is called {\em homogeneous homomorphism}, in short {\em h-homomorphism}. 
}\end{Def}
Analogous terminology will be used for the corresponding
Lie algebra homomorphisms of graded algebras that commute
with dilations.
\begin{Def}[P-differentiability]\label{Pansudiff}{\rm
Let $d$ and $\rho$ be homogeneous distances of $\G$ and $\M$,
respectively. We consider the mapping $f:\Omega\lra\M$.
We say that $f$ is {\em P-differentiable} at $x\in\Omega$ if
there exists an h-homomorphism $L:\G\lra\M$ such that
$$
\frac{\rho\big(f(x)^{-1}f(xh),L(h)\big)}{d(h)}\lra0
\quad\mbox{as}\quad h\ra e\,.
$$
The h-homomorphism $L$ satisfying this limit is unique
and it is called {\em P-differential} of $f$ at $x$.
We denote $L$ by $Df(x)$, when we read the P-differential
between the corresponding Lie algebras, we will denote it by $df(x)$.
}\end{Def}
\begin{The}\label{Pansu}
Every Lipschitz mapping $f:\Omega\lra\M$
is $\mu$-a.e. P-differentiable.
\end{The}
This theorem is an important result due to Pansu, \cite{Pan2}. 
Here have presented a slightly more general version where
$\M$ is graded, but it might not be stratified, \cite{Mag11}.
\begin{Def}[Distributional derivatives]\label{distder}{\rm
Let $\Omega$ be an open subset of a stratified group $\G$, let $X$ be a 
left invariant vector field of $\G$ and let $E$ be a
finite dimensional normed space.
Then for every $F\in L^1_{loc}(\Omega,E)$ we say that 
$G\in L^1_{loc}(\Omega,E)$ is the {\em distributional derivative}
of $F$ with respect to $X$ if 
\[
\int_\Omega F\; X\ph\;d\mu =-\int_\Omega G\;\ph\;d\mu
\]
for every $\ph\in C_c^\infty(\Omega)$. 
Uniqueness of $G$ allows us to use the notation $\nabla_XF$.
In the case $\Omega$ is an open subset of $\R^n$ and
$X=\der_{x_j}$ we will use the notation $\nabla_j$,
where $e_j$ belongs to the canonical basis of $\R^n$.
}\end{Def}

\begin{Rem}{\rm
To avoid confusion, we stress that the symbol $\nabla_X$
will always denote a distributional derivative, since we
will never consider connections in this work.}
\end{Rem}
\begin{Def}[Horizontal gradient]{\rm
Under the conditions of Definition~\ref{distder}
if we have equipped $V_1$ with a scalar product
and $(X_1,\ldots,X_m)$ is an orthonormal basis of $V_1$,
then we introduce the notation
\[
\nabla_HF=\big(\nabla_{X_1}F,\ldots,\nabla_{X_m}F\big)
\]
to denote the {\em distributional horizontal gradient} of $F:\Omega\lra E$.
}\end{Def}
\subsection{Some auxiliary results}
In this subsection $\M$ denotes a graded group equipped with Lie algebra $\cM$.
The next lemma easily follows from the scaling property of the homogeneous distance
in a graded group.
\begin{Lem}
Let $\cM=W_1\oplus\cdots\oplus W_\upsilon$ and
let $\pi^i:\cM\lra W_i\oplus\cdots\oplus W_\upsilon$
be the canonical projection. Let $U$ be a bounded open neighbourhood
of the unit element $e\in\M$.
Then there exists a constant $K_U>0$, depending on $U$, such that
\begin{eqnarray}\label{estimil}
\|\pi^i\left(\exp^{-1}(x)\right)\|\leq K_U\;d(x)^i	
\end{eqnarray}
holds for every $x\in U$ and every $i=1,\ldots,\upsilon$.
\end{Lem}
\begin{Rem}[Contact property]\label{cntpro}{\rm
We wish to point out how the contact property of a mapping
$f:\Omega\lra\M$ at some point $x\in\Omega$, namely,
\[
df(x)(H_x\G)\subset H_{f(x)}\M
\]
is equivalent to the differential constraint
\begin{eqnarray}\label{xifj}
X_iF_j-\sum_{n=2}^\upsilon\frac{(-1)^n}{n!}\,
\pi_j\big([F_{j-1}(x),X_iF_{j-1}(x)]_{n-1}\big)=0\,,
\end{eqnarray}
for every $i=1,\ldots,m$ and every $j=2,\ldots,\up$,
where $(X_1,\ldots,X_m)$ is a basis of $V_1$, $f=\exp\circ F$
and $F_j=\pi_j\circ F$. 
The exponential mapping $\exp:\cM\lra\M$ satisfies
\begin{equation}\label{dexp}
d\exp\,(X)=\mbox{Id}-\sum_{n=2}^\upsilon\frac{(-1)^n}{n!}\,\mbox{ad}(X)^{n-1}\,,
\end{equation}
see Theorem~2.14.3 of \cite{Vara}.
Identifying for every $y\in\M$ the tangent space $T_y\M$ with $\cM$
and applying formula \eqref{dexp}, we have
\begin{eqnarray*}
X_iF(x)=X_iF(x)-\sum_{n=2}^\upsilon\frac{(-1)^n}{n!}\,
\mbox{ad}\big(F(x)\big)^{n-1}\big(X_iF(x)\big)\,.
\end{eqnarray*}
Then $X_iF(x)\in H_{f(x)}\M$ if and only if
$$
\pi_j\Big(X_iF(x)-\sum_{n=2}^\upsilon\frac{(-1)^n}{n!}
[F(x),X_iF(x)]_{n-1}\Big)=0\quad\mbox{for all}\quad j\geq2\,.
$$
This proves the characterizing property of equations \eqref{xifj}.
We also notice that this characterization holds even if the 
mapping is differentiable along horizontal directions,
see \cite{Mag11} for the precise definition of horizontal differentiability.
}\end{Rem}
The following result corresponds to Corollary~5.4 of \cite{Mag11}.
\begin{The}\label{hdiffcurve}
Let $\Gamma:[a,b]\lra\M$ be a curve and define $\gamma=\exp^{-1}\circ\Gamma=\sum_{i=1}^\upsilon\gamma_i$,
where $\gamma_i$ takes values in $W_i$.
Then the following statements are equivalent:
\begin{enumerate}
	\item $\Gamma$ is Lipschitz continuous,
	\item $\gamma$ is Lipschitz continuous and
	the differential equation
	\begin{eqnarray}\label{eqhdiffcurve1}
	\dot\gamma_i(t)=\sum_{n=2}^\upsilon\frac{(-1)^n}{n!}\;
	\pi_i\left([\gamma(t),\dot\gamma(t)]_{n-1}\right)
	\end{eqnarray}
	is a.e. satisfied for every $i\geq2$.
\end{enumerate}
If one of the previous conditions holds, then
there exists a constant $C>0$ only depending
on $\rho$ and $\|\cdot\|$,
such that for any $\tau_1<\tau_2$, we have
\begin{eqnarray}\label{LipEstCurve}
\rho\big(\Gamma(\tau_1),\Gamma(\tau_2)\big)\leq 
C\,\int_{\tau_1}^{\tau_2}\|\dot\gamma_1(t)\|\,dt\,.
\end{eqnarray}
\end{The}
\begin{Rem}{\rm
If $\M$ is a stratified group, namely, its horizontal distribution
satisfies the Lie bracket generating condition, then the previous
theorem is well known.
In fact, in this case equivalence of (1) and (2) would
follow by Proposition~11.4 of \cite{HK} joined with Remark~\ref{cntpro}.
This proposition shows that Lipschitz curves can be 
characterized as horizontal curves in the more general
Carnot-Carath\'eodory spaces, that clearly include stratified groups.
Moreover, estimate \eqref{LipEstCurve} can be obtained for
instance by the sub-Riemannian area formula \cite{Mag}
joined with Theorem~2.10.13 of \cite{Fed}.
However, graded groups need not satisfy the Lie bracket generating condition,
then they are not included in the family of Carnot-Carath\'eodory spaces.
}\end{Rem} 
%
%
%
%%%%%%%%%%%%%%%%%%%%%%%%%%%%%%%%%%%%%%%%%%%%%%%%%%%%%%%%%%%%%%%%%%%%%%
%
%
%
%
\section{Technical lemmata}
%
%
%
%
%%%%%%%%%%%%%%%%%%%%%%%%%%%%%%%%%%%%%%%%%%%%%%%%%%%%%%%%%%%%%%%%%%%%%%%
%
%
%
In this section we study the properties of mappings
$f:A\lra\M$ where $\M$ is a graded group and $A$ may vary.
$F$ denotes the mapping $\exp^{-1}\circ f:A\lra\cM$,
where $\cM$ is the Lie algebra of $\M$.
The canonical projections onto the layers $W_i$ 
of the algebra $\cM$ are the mappings $\pi:\cM\lra\ W_i$.
We use the notation
$F_i=\pi\circ F$, hence
\[
F=F_1+\cdots+ F_\up\quad \mbox{and} \quad
f=\exp\big(F_1+\cdots +F_\up\big)\,.
\]

\begin{Lem}\label{uX}
Let $X\in V_1$ and let $N\subset\G$ be a normal subgroup such that
$\G=N\rtimes H_X$.
Let $O\subset N$ and $J\subset H_X$ be open subsets, where $J$ is connected,
and consider $u\in L^1_{loc}(\Omega)$, where $\Omega=OJ$ is an open set.
If the distributional derivative $D_Xu$ belongs to $L^\infty(\Omega)$,
then up to redefinition of $u$ an a $\mu$-negligilbe set,
for $\mu_N$-a.e. $n\in O$, we have 
\[
|u\big(n\exp(tX)\big)-u\big(n\exp(\tau X)\big)|\leq \|D_Xu\|_{L^\infty(\Omega)}
\;|t-\tau|\,.
\]
\end{Lem}
{\sc Proof.}
We adopt the same notation used in the proof of Proposition~\ref{FubX}.
Then we recall tha mapping $\psi:\R\times\R^{q-1}\lra\M$,
observing that $\tilde u=u\circ\psi$ belongs to $L^1_{loc}\big(\psi^{-1}(\Omega)\big)$.
Taking into account that $F^{-1}\circ\psi$ has jacobian equal to one,
by definition of distirbutional derivative one easily gets the following
equality of distributional derivatives
\[
\nabla_1\tilde u=(D_Xu)\circ\psi\,.
\]
Then $\tilde u\in L^1_{loc}\big(\psi^{-1}(\Omega)\big)$
and $\nabla_1\tilde u\in L^\infty\big(\psi^{-1}(\Omega)\big)$.
The set $\psi^{-1}(O)=\widetilde O$ is an open subset of $\{0\}\times \R^{q-1}$
and $\widetilde J=\psi^{-1}(J)$ is an open interval of $\R\times\{0\}\subset\R^q$,
hence $\psi^{-1}(O)=\widetilde J\times\widetilde O$.
Thus, by a standard mollification argument, see for instance \cite{EvansGar}, Theorem 2 of Section 4.9.2, our claim follows. $\Box$
\begin{Lem}\label{OJ}
Let $X\in V_1$ and let $\G=N\rtimes H_X$, where $N$ is a homogeneous normal subgroup.
Let $O$ and $J$ be open subsets of $N$ and of $H_X$, respectively, where $J$ is
connected, and let $z\in\G$.
We consider the open set $\Omega=zOJ$ along with the continuous mapping $f:\Omega\lra\M$. Let $(X_1,\ldots,X_m)$ be a basis of $V_1$ and assume that there exist
\begin{eqnarray}\label{D_XinLinfty}
\nabla_{X_i}F_j\in L^\infty_{loc}(\Omega,W_j)\quad\mbox{and}\quad
\nabla_{X_i}F_1\in L^\infty(\Omega,W_1)
\end{eqnarray}
for every $i=1,\ldots,m$ and $j=2,\ldots,\upsilon$ and the contact equations
\begin{equation}\label{XF_j}
X_iF_j=\sum_{n=2}^\upsilon\frac{(-1)^n}{n!}\;\pi_j\left([F,X_iF]_{n-1}\right)\\
\end{equation}
$\mu$-a.e. hold in $\Omega$, then
there exists $C>0$, depending on $\rho$ and $\|\cdot\|$, such that
\[
\rho\lls f\big(zn\exp(tX)\big),f\big(zn\exp(\tau X)\big)\rls\leq C\,
\|XF_1\|_{L^\infty(\Omega)} |t-\tau|
\]
for every $\exp(tX),\exp(\tau X)\in J$ and every $n\in O$.
\end{Lem}
{\sc Proof.}
Up to a left translation, taking into account the left invariance
of $X$, it is not restrictive to assume that $z$ is the unit element.
Then we consider the curve
\[
\Gamma_{n}(s)=n\exp(sX)=\exp\big(\gamma_{n,1}(s)+\cdots\gamma_{n,\upsilon}(s)\big)
\] 
where $\gamma_{n,j}=\pi_j\circ\gamma_n$ and $\gamma_n=\exp^{-1}\circ \Gamma_n$.
We choose $t,\tau\in\R$ with $t<\tau$ such that $\exp(tX),\exp(\tau X)\in J$.
By linearity, we have $D_XF$ exists and belongs to $L^\infty_{loc}(\Omega,\cM)$.
Taking into account continuity of all $F_j$'s, Lemma~\ref{uX} gives
\[
\|F_j\big(\Gamma_n(t)\big)-F_j\big(\Gamma_n(\tau)\big)\|\leq c\,\|XF_j\|_{L^\infty(\Omega'_n)}
\,(\tau-t)
\]
for all $n\in O$, where $c>0$ depends on $\|\cdot\|$ and
$\Omega'_n$ is an open neighbourhood of $\Gamma_n([t,\tau])$, that is
compactly contained in $\Omega$.
In particular, $t\ra F\circ\Gamma_n(t)$ is absolutely continuous 
on compact intervals for all $n\in O$.
By Proposition~\ref{FubX}, from the $\mu$-a.e. validity of \eqref{XF_j} it follows that
for $\mu_N$-a.e. $n\in O$, we have
\[
\frac{d}{ds}\Big(F_j\circ\Gamma_n\Big)(s)=XF_j\circ\Gamma_n(s)=
\sum_{n=2}^\upsilon\frac{(-1)^n}{n!}\;\pi_j
\left([F\circ\Gamma_n(s),(XF)\circ\Gamma_n(s)]_{n-1}\right)
\]
for a.e. $s\in\{l\in\R\mid \exp(lX)\in J\}$.
Thus, we can apply Theorem~\ref{hdiffcurve} to the curve
\[
f\circ \Gamma_n=\exp\big(F_1\circ\Gamma_n+\cdots+F_\upsilon\circ\Gamma_n\big)\,,
\]
getting
\begin{eqnarray}
\rho\big(f\circ\Gamma_n(t),f\circ\Gamma_n(\tau)\big)\leq C\,\int_{t}^{\tau}\|(XF_1)
\circ\Gamma_n(s)\|\,ds\,.
\end{eqnarray}
Then the continuity of $f$ and the hypothesis $XF_1\in L^\infty(\Omega)$
lead us to the conclusion. $\Box$
\begin{Lem}\label{LipF_1}
Let $f:X\lra\M$ be a Lipschitz mapping, where $X$ is a metric space.
Then there exists a constant $C>0$, depending on the norm $\|\cdot\|$ of $\cM$ and the distance $\rho$ of $\M$ such that
\[
\|F_1(x)-F_1(y)\|\leq C\,\Lip(f)\,d(x,y) \quad\mbox{for every}\quad x,y\in X\,.
\]
\end{Lem}
{\sc Proof.}
We define 
\begin{eqnarray}\label{homnorm}
|\xi|=\sum_{j=0}^\upsilon\|\xi_j\|^{1/j}
\end{eqnarray}
where $\xi=\xi_1+\xi_2+\cdots+\xi_\iota$ and $\xi_j\in W_j$.
By continuity and homogeneity one can check that there exists a constant $C>0$
such that
\[
\frac{1}{C}\,\rho\big(\exp\xi\big)\leq |\xi|\leq C\,\rho\big(\exp\xi\big)
\]
Then we have
\[
\|\xi_1-\eta_1\|\leq |-\xi\gope\eta|\quad\mbox{for every}\quad
\xi,\eta\in\cM,
\]
where $\xi_1=\pi_1(\xi)$ and $\eta_1=\pi_1(\eta)$.
It follows that
\begin{eqnarray}
	\|F_1(x)-F_1(y)\|\leq C\,\rho\big(f(x),f(y)\big)\leq
	C\,\Lip(f)\,d(x,y)\,.
\end{eqnarray}
This concludes the proof. $\Box$
\begin{Def}[Piecewise horizontal line]{\rm
A continuous curve $\Gamma:[a,b]\lra\G$ is a {\em piecewise horizontal line}
if there exist $n\in\N$, numbers $a\leq t_0<\cdots<t_n\leq b$
and $X_k\in V_1$ with $k=1,\ldots,n$ such that $\Gamma_{|[t_{k-1},t_k]}(t)=\Gamma(t_{k-1})\exp\big((t-t_{k-1})X_k\big)$.}
\end{Def}
\begin{Lem}\label{HorLinApp}
Let $\Omega$ be an open subset of $\G$ and let $\Gamma:[a,b]\lra\Omega$
be a Lipschitz curve. Then there exists a sequence of piecewise horizontal lines 
uniformly converging to $\Gamma$, whose lengths are also converging
to the length of $\Gamma$.
\end{Lem}
{\sc Proof.}
By Lemma~\ref{LipF_1}, it follows that the curve
$\gamma_1=\pi_1\circ \exp^{-1}\circ \Gamma:[a,b]\lra W_1$
is Lipschitz.
The derivative $\dot\gamma_1$ is essentially bounded.
Let $\ph_k:[a,b]\lra W_1$ be a sequence of piecewise constant functions 
that a.e. converge to $\dot\gamma_1$ and
$|\ph_k(t)|\leq|\dot\gamma_1(t)|$ for every $k$ and a.e. $t\in[a,b]$.
We define the Lipschitz functions
\[
\gamma_{1k}(t)=\gamma_1(0)+\int_0^t\ph_k(s)\,ds
\]
that uniformly converge to $\gamma_1$. 
Exploiting \eqref{eqhdiffcurve1}, we can define
\[
\gamma_{2k}(t)=\gamma_2(0)+\frac{1}{2}
\int_0^t\left[\gamma_{1k}(s),\dot\gamma_{1k}(s)\right]\,ds
\]
and more generally by iteration
\[
\gamma_{ik}(t)=\gamma_{ik}(0)+\sum_{n=2}^\upsilon\frac{(-1)^n}{n!}\;
\int_0^t\pi_i\bigg(\bigg[\sum_{j=1}^{i-1}\gamma_{jk}(s),\sum_{j=1}^{i-1}\dot\gamma_{jk}(s)\bigg]_{n-1}\bigg)\,ds.
\]
The curve $\Gamma_k=\exp\big(\sum_{j=1}^\upsilon\gamma_{jk}\big)$
is the unique horizontal lifting of $\gamma_k$.
We wish to show that it is piecewise horizontal and uniformly converge to $\Gamma$.
To see this, we observe that the projection on the first layer
of any curve
\[
t\lra \exp\xi\exp tX
\]
with $X\in V_1$ and $\xi\in\cG$ has the form $t\lra \xi_1+tX$, where $\pi_1(\xi)=\xi_1\in V_1$. We notice that $\gamma_{1k}$ has exactly this
form, since $\ph_k$ is piecewise constant, then the uniqueness of the horizontal lifting implies that $\Gamma_k$ is piecewise horizontal.
By construction, the fact that $\gamma_{1k}$ uniformly converges to $\gamma_1$
implies the uniform convergence of $\Gamma_k$ to $\Gamma$.
From Corollary~5.5 of \cite{Mag11}, we get the formula
\[
\mbox{\rm length}(\Gamma_k)=
\mbox{\rm Var}_a^b\Gamma=\int_a^b\rho\Big(\exp\big(\dot\gamma_{1k}(t)\big)\Big)\,dt\,.
\]
It follows that $\mbox{\rm length}(\Gamma_k)\ra\mbox{\rm length}(\Gamma)$
as $k\ra\infty$. $\Box$
%
%
%
%%%%%%%%%%%%%%%%%%%%%%%%%%%%%%%%%%%%%%%%%%%%%%%%%%%%%%%%%%%%%%%%%%%%%%
%
%
%
%
\section{Group-valued mappings on Riemannian manifolds}
%
%
%
%%%%%%%%%%%%%%%%%%%%%%%%%%%%%%%%%%%%%%%%%%%%%%%%%%%%%%%%%%%%%%%%%%%%%%
%
%
%
%

In this section we extend contact equations to graded group-valued mappings on
Riemannian manifolds and find the corresponding differential characterization of
the Lipschitz property.
In this case, contact equations can be written using the standard differential
for differentiable manifolds.

Throughout this section, we have the following assumptions.
We denote by $\M$ a graded group with graded algebra $\cM$
and $N$ indicates a Riemannian manifold.
The symbol $\rho$ denotes a homogeneous distance of $\M$.
Any function $f: N\lra\M$ is also written as the composition $f=\exp\circ F$,
where $F:N\lra\cM$ and $\exp:\cM\lra\M$ denotes the exponential of Lie groups.
We will also use the notation $F_i=\pi\circ F$, where $F=F_1+\cdots+ F_\up$,
$f=\exp\big(F_1+\cdots +F_\up\big)$ and $\pi_i:\cM\lra W_i$ is the canonical
projection onto the $i$-th layer.
\begin{Pro}\label{1DLipContact}
Let $N$ be a one dimensional connected Riemannian manifold
and let $F:N\lra\cM$ be a Lipschitz mapping such that
\begin{equation}\label{1Drkcont}
\dot F_j=\sum_{n=2}^\upsilon\frac{(-1)^n}{n!}\;\pi_j\left([F,\dot F]_{n-1}\right)\\
\end{equation}
a.e. holds in $N$ for every $j=2,\ldots,\upsilon$.
Then the associated mapping $f:N\lra\M$ is Lipschitz continuous and 
there exists a geometric constant $C>0$
depending on the norm of $\cM$ and the distance of $\M$ such that
$ \Lip(f)\leq C\;\Lip(F_1). $
\end{Pro}
{\sc Proof.}
Let $x,y$ be two points of $N$ and let $c:[0,L]\lra N$ be
the length minimizing geodesic such that $c(0)=x$, $c(L)=y$ and $L=d(x,y)$.
The curve $\gamma=F\circ c$ is Lipschitz and in view of \eqref{1Drkcont}
it satisfies \eqref{eqhdiffcurve1} for a.e. $t$ in $[0,L]$, then
by Theorem~\ref{hdiffcurve}, we get
\begin{equation}\label{rhof}
\rho\big(f(x),f(y)\big)=\rho\big(\Gamma(0),\Gamma(L)\big)\leq
C\,\int_0^L\|(F_1\circ c)'(t)\|\,dt
\end{equation}
for a suitable geometric constant $C>0$ depending on $\M$ and $\cM$,
where we have set
$\Gamma=f\circ c$. Since
\[
\big|(F_1\circ c)(t+h)-(F_1\circ c)(t)\big|\leq\Lip(F_1)\,d\big(c(t+h),c(t)\big)=
\Lip(F_1)\,|h|\,,
\]
where $d$ is the Riemannian distance of $N$.
Then \eqref{rhof} leads us to our claim. $\Box$
\begin{Rem}{\rm
As the linearity of $X\lra[Z,X]_{n-1}=\mbox{ad}^{n-1}(X)$ makes 
\eqref{1Drkcont} intrinsic, the same holds for higher dimensional manifolds.
In fact, if $N$ has dimension $k\geq1$ and $f$ is differentiable at a point $\overline{x}$ of $N$, then the contact condition $df(\overline{x})(T_{\overline{x}}N)\subset H_{f(\overline{x})}\M$
is equivalent to the validity of the system
\begin{equation}\label{rmce}
\partial_{x_l}F_j(\overline{x})=\sum_{n=2}^\upsilon\frac{(-1)^n}{n!}\;\pi_j
\left([F(\overline{x}),\partial_{x_l}F(\overline{x})]_{n-1}\right)\,,
\end{equation}
for every $j=2,\ldots,\upsilon$ and $l=1,\ldots,k$. 
This follows by Remark~\ref{cntpro}. Since the previous formula
is independent of the local coordinates $(x_l)$ chosen around $\overline{x}$ in $N$,
then an equivalent intrinsic version is the following one
\begin{equation}\label{rmce1}
dF_j(\overline{x})=\sum_{n=2}^\upsilon\frac{(-1)^n}{n!}\;\pi_j
\left([F(\overline{x}),dF(\overline{x})]_{n-1}\right)\,,
\end{equation}
for every $j=2,\ldots,\upsilon$. Here $d$ denotes the standard differential
for mappings on differentiable manifolds and
\[
dF_j(\overline{x}):T_{\overline{x}}N\lra H^j_{f(\overline{x})}\M.
\]
Formula \eqref{rmce} can be also written in a shorte and intrinsic form
adopting the Riemannian gradient as follows
\begin{equation}
\nabla F_j(\overline{x})=\sum_{n=2}^\upsilon\frac{(-1)^n}{n!}\;\pi_j
\left([F(\overline{x}),\nabla F(\overline{x})]_{n-1}\right)\,.
\end{equation}
}\end{Rem}
\begin{Def}{\rm
Let $N$ be a Riemannian manifold and let $d$ denote the Riemannian distance.
We say that a subset $O\subset N$ is {\em geodetically convex} if for every
$x,y\in O$, there exists a length minimizing geodesic $c:[0,L]\lra O$
such that $c(0)=x$, $c(L)=y$.
}\end{Def}
\begin{Rem}{\rm
Notice that notion of geodetic convexity we are adopting is not standard,
since in Riemannian Geometry different forms of uniqueness of the 
connecting and minimizing geodesic are also required, see for instance \cite{Chavel}.
For instance, according to our definition, the intersection of the
$2$-dimensional sphere $S^2$ embedded in $\R^3$ with any closed
half space, where $S^2$ is equipped with the canonical Riemannian metric,
is an example of geodetically convex set.}
\end{Rem}
\begin{The}\label{LipContact}
Let $N$ be a connected Riemannian manifold of dimension
higher than one and let $O\subset N$ be an open and geodetically convex set.
Let $\nu$ be the canonical Riemannian measure on $N$ and let
$F:O\lra\cM$ be a Lipschitz mapping such that
\begin{equation}\label{rkcont}
dF_j=\sum_{n=2}^\upsilon\frac{(-1)^n}{n!}\;\pi_j\left([F,dF]_{n-1}\right)\\
\end{equation}
$\nu$-a.e. holds in $O$ for every $j=2,\ldots,\upsilon$.
Then the associated mapping $f:O\lra\M$ is Lipschitz continuous and 
there exists a geometric constant $C>0$
depending on the norm of $\cM$ and the distance of $\M$ such that
$ \Lip(f)\leq C\;\Lip(F_1). $
\end{The}
{\sc Proof.}
We choose two arbitrary points $p,q\in O$ and consider the smooth length minimizing 
geodesic $c:[0,L]\lra N$, where $L=d(p,q)$, $c(0)=p$ and $c(L)=q$. Here $d$ denotes the
Riemannian distance on $N$. 
Let $\tau\in[0,L]$ and let $\big(c(\tau),\dot c(\tau)\big)\in TN$.
Let $\cT N$ be the open subset of $TN$ corresponding to the domain of the
Riemannian exponential mapping $\cExp:\cT N\lra N$, see Theorem~I.3.2 of \cite{Chavel}.
We wish to construct a tubular neighbourhood of geodesics containing also the image of
the restriction of $c$ to a neighbourhood of $\tau$.

To do this, we fix $z=c(\tau)$, $\xi=\dot c(\tau)$
and select an arbitrary embedded smooth one codimensional
submanifold $\Sigma\subset N$ passing through $z$, such that $T_z\Sigma$
is orthogonal to $\xi\in T_zN$.
We will select all geodesics tangent to a normal
field of $\Sigma$ in a neighbourhood of $z$. To make this argument rigorous,
according to Chapter 5, p.132, of \cite{Petersen},
we consider the vector bundle $T\Sigma^\bot$ of fibers
\[
T_s\Sigma^\bot=\big\{v\in T_sN\mid v\in(T_s\Sigma)^\bot\subset T_sN\big\}
\quad \mbox{for every $s\in\Sigma$.}
\]
By definition, we have the orthogonal decomposition
$T_sN=T_s\Sigma\oplus T_s\Sigma^\bot$.
Now, we consider the normal exponential mapping
$\cExp^\bot$ as the following restriction
\[
\cExp^\bot:\cT N\cap T\Sigma^\bot\lra N\,.
\]
Since the differential of $\cExp^\bot$ is nonsingular at every point $(s,0)$,
in particular there exist open neighbourhoods $\cU$ of $(z,0)$ in $\cT N\cap T\Sigma^\bot$
and $U$ of $z$ in $N$ such that $\cExp^\bot_{|}:\cU\lra U$ is a smooth diffeomorphism.
This provides us with a local system of coordinates around $z$ made by the
local geodesic flow. Let $k$ be the dimension of $N$, hence up to shrinking
both $\cU$ and $U$, we can select local coordinates $(y_1,\ldots,y_{k-1})$ of 
$\Sigma$ centered at zero, around $z$ and fix the local unit normal field $n$ of 
$\Sigma$ around $z$ such that ${\bf n}(z)=\xi$.
Then we define 
\[
H(y,t)=\cExp^\bot_{|}\Big(\zeta(y),t\,{\bf n}\big(\zeta(y)\big)\Big)
\]
where $(A,\zeta)$ is a local chart of $\Sigma$,
$A$ is an open subset of $\R^{k-1}$ containing the origin, with $\zeta(0)=z$, and
$\Big(\zeta(y),t\,{\bf n}\big(\zeta(y)\big)\Big)\in\cU$ if and only if
$(y,t)\in A\times I$ for a suitable open interval $I$ of $\R$.
By local uniqueness of geodesics, we get
\[
I\ni t\lra c(\tau+t)=H(0,t)=\cExp^\bot\big(z,t\,{\bf n}(z)\big)
=\cExp\big(z,t\,\xi\big)\,.
\]
Since $H:A\times I\lra U$ is bilipschitz, then in view of 
Fubini's theorem for a.e. $y\in A$ we have that 
$F$ is differentiable at $H(y,t)$ for a.e. $t\in I$ and there satisfies
\eqref{rkcont} for every $j=2,\ldots,\upsilon$.
Then in particular, we have the partial derivatives
\begin{equation}
\der_t(F_j\circ H)(y,\cdot)
=\sum_{n=2}^\upsilon\frac{(-1)^n}{n!}\;\pi_j
\left(\big[(F\circ H)(y,\cdot),\der_t(F_j\circ H)(y,\cdot)\big]_{n-1}\right)\\
\end{equation}
a.e. in $I$ for a.e. $y\in A$.
Since $F\circ H(y,\cdot)$ is also Lipschitz, then for a.e. $y\in A$
Proposition~\ref{1DLipContact}
yields a geometric constant $C_1>0$ such that
\[
\rho\big(f\circ H(y,t),f\circ H(y,t')\big)\leq 
\Lip\big(F_1\circ H(y,\cdot)\big)\, C \,|t-t'|\,.
\]
Continuity extends the previous estimate to all $y\in A$ and $t\in I$.
In particular, for $y=0$,
 it follows that
\begin{equation}\label{liptau}
\rho\big((f\circ c)(\tau+t),(f\circ c)(\tau+t')\big)\leq 
C\, \Lip(F_1) \,|t-t'|
\end{equation}
for every $t,t'\in I$. The arbitrary choice of $\tau\in[0,L]$ 
gives a finite open covering of $[0,L]$ made of intervals satisfying \eqref{liptau}.
Since the constants of this estimate are independent of $\tau$, this
leads us to the end of the proof. $\Box$
\begin{Rem}\label{LipContactCl}{\rm
Under the hypotheses of the previous theorem, if we assume in addition
that $N$ is a complete Riemannian manifold, then the mapping $f$ extends 
to a Lipschitz mapping $f:\overline{O}\lra\M$, 
with $\Lip(f)\leq C\;\Lip(F_1).$ Since compact Riemannian manifolds are complete
and then geodetically convex, then Theorem~\ref{LipContact}
obviously implies Theorem~\ref{RMLip}.
}\end{Rem}
%
%
%
%
%
%
%
%%%%%%%%%%%%%%%%%%%%%%%%%%%%%%%%%%%%%%%%%%%%%%%%%%%%%%%%%%%%%%%%%%%%%%
%
%
%
%
\section{The differential characterization}
%
%
%
%
%%%%%%%%%%%%%%%%%%%%%%%%%%%%%%%%%%%%%%%%%%%%%%%%%%%%%%%%%%%%%%%%%%%%%%%
%
%
%
%
%
%
%

In this section we give a proof of the differential characterization
of locally Lipschitz mappings and show a simple application.
\vskip.2cm
{\sc Proof of Theorem~\ref{CharLip}.}
Let $f$ be locally Lipschitz.
By Lemma~\ref{LipF_1}, there exists $C>0$ only depending on the norm
$\|\cdot\|$ and the distance $\rho$, such that
\begin{eqnarray}
	\|F_1(x)-F_1(y)\|\leq C\,\rho\big(f(x),f(y)\big)
\end{eqnarray}
for every $x,y\in\Omega$. Then $F_1$ is also locally Lipschitz.
We wish to show that all $F_j$'s are locally Lipschitz. 
By the Baker-Campbell-Hausdorff formula \eqref{absBCH}, we have
\[
\|-\xi_i+\eta_i+\sum_{n=2}^\upsilon\pi_i\big(c_n(-\xi,\eta)\big)\|=\left|\pi_i\big(-\xi\gope\eta\big)\right|^i\leq
\Big[C\;\rho(\exp\xi,\exp\eta)\Big]^i,
\]
then estimate \eqref{bilestim} yields $\tilde\alpha_n(\nu)>0$,
with $\nu=\max\{\|\xi\|,\|\eta\|\}$, such that
\[
\|-\xi_i+\eta_i\|\leq\sum_{n=2}^\upsilon\tilde\alpha_n(\nu)
\|[\xi,\eta]\|+\Big[C\;\rho(\exp\xi,\exp\eta)\Big]^i,
\]
then taking into account \eqref{lcst} and \eqref{estimil}, we get
\[
\|-\xi_i+\eta_i\|\leq\left(\sum_{n=2}^\upsilon\tilde\alpha_n(\nu)\,
\beta\,\nu\, K_U
+C^i\,\rho(\exp\xi,\exp\eta)^{i-1}\right)\rho(\exp\xi,\exp\eta)
\]
where $U$ is a compact set containing $\exp\xi$ and $\exp\eta$,
hence depending on $\nu$.
As a consequence, replacing both $\xi_i$ and $\eta_i$ with $F_i(x)$ and $F_i(y)$,
we have shown that all $F_i$'s are locally Lipschitz.
As a standard fact, it follows that all $\nabla_{X_i}F_j$ belong to $L^\infty_{loc}(\Omega)$.
To prove the a.e. validity of \eqref{hcondF}, we choose an arbitrary
left invariant vector field $X\in V_1$ and by Lemma~\ref{factorX}
we get a normal subgroup $N$ of $\G$ such that $\G=N\rtimes H_X$.
We can cover $\Omega$ by a countable union of open sets of the form $zOJ$, where
$z\in\G$, $O$ and $J$ are connected open neighbourhoods of the unit
element in $N$ and in $H_X$, respectively.
Thus, we can reduce ourselves to prove our claim in the case $\Omega=zOJ$.
Due to Theorem~\ref{Pansu} in the case $\M$ is a vector space,
we consider the full measure set $A\subset\Omega$ of points where all $F_j$'s are P-differentiable. Then formula \eqref{FubNX} implies
\[
\nu_X\big(J\sm p^{-1}z^{-1}A\big)=0
\]
for $\mu_N$-a.e. $p\in O$.
Let us pick one of these $p$'s. For a.e. $t$ we have
$zp\exp(tX)\in A$, then all $F_j$'s are P-differentiable
at this point. Definition of P-differentiability yields
\[
DF_j\big(c_{zp,X}(t)\big)(X)=XF_j\big(c_{zp,X}(t)\big),
\]
for every $j=2,\ldots,\upsilon$, where we have defined
\[
\tau\lra c_{zp,X}(\tau)=zp\,\exp\big(\tau X\big).
\]
As a result, setting
\[
\tau\lra\Gamma_{zp,X}(\tau)=f\circ c_{zp,X}(\tau)
=\exp\gamma_{zp,X}(\tau)=\exp\sum_{j=1}^\upsilon \gamma_{zp,X,j}(\tau),
\]
where we have set $\gamma_{zp,X,j}=\pi_j\circ \gamma_{zp,X}$.
We have then proved that 
\begin{eqnarray}\label{Xfj}
\dot\gamma_{zp,X,j}(t)=XF_j\circ c_{zp,X}(t).
\end{eqnarray}
By Theorem~\ref{hdiffcurve}, the local Lipschitz property of $f$ implies 
the a.e. validity of \eqref{eqhdiffcurve1} for the curve $\gamma_{zp,X}$.
Thus, taking into account \eqref{Xfj}, we obtain
\[
XF_j(c_{zp,X}(\tau))=\sum_{n=2}^\upsilon\frac{(-1)^n}{n!}\;
\pi_j\left([F(c_{zp,X}(\tau)),XF(c_{zp,X}(\tau))]_{n-1}\right),
\] 
for a.e. $\tau$.
We have proved the a.e. validity of the previous equation for
$\mu_N$-a.e. $p\in O$, therefore Proposition~\ref{FubX} implies the validity of 
\begin{eqnarray}\label{XFj}
XF_j=\sum_{n=2}^\upsilon\frac{(-1)^n}{n!}\;\pi_j\left([F,XF]_{n-1}\right)	
\end{eqnarray}
$\mu$-a.e. in $\Omega$ for every $j=2,\ldots,\upsilon$.
The arbitrary choice of $X$ gives the validity of \eqref{hcondF}.

Conversely, we assume the validity of the second condition.
By linearity of distributional derivative and of \eqref{hcondF}
for every $X\in V_1$, we have $D_XF\in L^\infty_{loc}(\Omega,\cM)$ and
\eqref{XFj} $\mu$-a.e. holds for every $j=2,\ldots,\upsilon$.
We fix $\|X\|=1$ and choose $z\in\Omega$.
We select connected open neighbourhoods $O\subset N$ and $J\subset H_X$
of the unit element $e\in\G$ such that
$\Omega'=zOJ$ is compactly contained in $\Omega$.
Thus, we can apply Lemma~\ref{OJ}, getting
\begin{eqnarray}\label{horestline}
&&\rho\lLs f\big(zn\exp(tX)\big),f\big(zn\exp(\tau X)\big)\rLs\leq
C'\,\|XF_1\|_{L^\infty(\Omega)}\; |t-\tau|\\
&&\leq C'\,\|\nabla_HF_1\|_{L^\infty(\Omega)}\; |t-\tau| \nonumber
\end{eqnarray}
for every $n\in O$ and every $\exp(tX),\exp(\tau X)\in J$.
Now, we consider the Carnot-Carath\'eodory distance
$\delta$ generated by a left invariant metric fixed on $\G$, see for instance \cite{Gr1}. This distance is homogeneous and then it is equivalent to $d$.
It is well known that $(\G,\delta)$ is a length space, namely 
any couple of points $x,y\in\G$ is connected by a geodesic whose length equals their distance $\delta(x,y)$.
Then from \eqref{horestline}, we get
\begin{equation}\label{LipLineEst}
	\rho\big(f(p),f(p')\big)
	\leq C\,\|\nabla_HF_1\|_{L^\infty(\Omega)}\,
	\delta(p,p')\qquad\mbox{for every}\qquad p,p'\in zJ\,,
\end{equation}
where $C=C'/\inf_{\|Y\|=1}\delta(\exp Y)$.
Now we arbitrarily choose $r>0$ and $p\in\Omega$ such that
$\overline{\cB_{p,6r}}\subset\Omega$,
where we have denoted by $\cB_{p,r}$ the open ball of center $x$ and 
radius $r$ with respect to $\delta$.
Let $x,y\in B_{p,2r}$ and let $\Gamma:[0,\delta(x,y)]\lra\G$
be the geodesic connecting $x$ with $y$. By triangle inequality,
it follows that the image of $\Gamma$ is contained in $B_{p,6r}$.
By Lemma~\ref{HorLinApp} we can find a sequence $(\Gamma_k)$ of piecewise horizontal lines defined in $[0,\delta(x,y)]$ and contained in $B_{p,6r}$
that uniformly converge to $\Gamma$ and their lengths converge to $\Gamma$.
On any horizontal segment of $\Gamma_k$ with horizontal direction
$X_i$, that is also a geodesic, we apply the estimate
\eqref{LipLineEst} where $J$ is considered contained in the subgroup
$H_{X_i}$. Thus, triangle inequality yields
\begin{eqnarray*}
\rho\lLs f\big(\Gamma_k(0)\big),f\big(\Gamma_k(\delta(x,y))\big)\rLs\leq C\,\|\nabla_HF_1\|_{L^\infty(\Omega)}\;l(\Gamma_k)\,,
\end{eqnarray*}
where $K_{p,r}=\overline{\cB_{p,6r}}$ and
$l(\Gamma_k)$ is the length of $\Gamma_k$ with respect to $\delta$.
Passing to the limit as $k\ra\infty$, we have shown that
\begin{eqnarray}
\rho\big(f(x)),f(y)\big)\leq C\,\|\nabla_HF_1\|_{L^\infty(\Omega)}\;\delta(x,y)
\end{eqnarray}
for every $x,y\in B_{p,2r}$. Adopting the same argument of Theorem~3.18
of \cite{Mag6}, it follows that $f$ is Lipschitz continuous on compact
sets of $\Omega$, namely, $f$ is locally Lipschitz.
Now, we show that the third condition is equivalent to the previous ones.
We first assume that the first condition holds, namely, $f$
is locally Lipschitz. Due to Theorem~\ref{Pansu}, $f$ is a.e.
P-differentiable and the equivalence of the first two conditions
clearly yields $\nabla_{X_i}F\in L^\infty_{loc}(\Omega,\cM)$
for every $i=1,\ldots,m$.
If we know that $f$ is a.e. P-differentiable and 
$\nabla_{X_i}F\in L^\infty_{loc}(\Omega,\cM)$
for every $i=1,\ldots,m$, then we apply 
Theorem~4.8 of \cite{Mag11}, according to which the pointwise
P-differentiability of $f$ implies the pointwise P-differentiability
of all $F_j:\Omega\lra\cM$ with the validity of formulae
\begin{eqnarray*}
&&\pi_1\circ df(x)=dF_1(x) \\
&&dF_i(x)(h)=\sum_{n=2}^\upsilon\frac{(-1)^n}{n!}\;\pi_i
\left([F(x),dF(x)(h)]_{n-1}\right)\,.
\end{eqnarray*}
at P-differentiability points $x$.
Taking into account that $dF_i$ denotes the P-differential
read in the Lie algebras and that $dF_i(x)(\exp X)=XF_i$,
these formulae implies the a.e. validity of contact equations
\eqref{hcondF}. We have then shown that that the third conditions
implies the second one. This concludes the proof. $\Box$
%
%
%
%%%%%%%%%%%%%%%%%%%%%%%%%%%%%%%%%%%%%%%%%%%%%%%%%%%%%%%%%%%%%%%%%%%%%%%%%%%%%%%%%%%%%%%
%
%
%
%
%
%
\subsection{Smooth functions that are not locally Lipschitz}\label{hNnLip}
%
%
%
%
%
%
%%%%%%%%%%%%%%%%%%%%%%%%%%%%%%%%%%%%%%%%%%%%%%%%%%%%%%%%%%%%%%%%%%%%%%%%%%%%%%%%%%%%%%%
%
%
%
We consider exponential coordinates $(x_1,x_2,x_3)$ of the Heisenberg group $\H^1$, 
with $p=\exp\big(\sum_{j=1}^3x_jX_j\big)\in\H^1$ and the basis of left invariant
vector fields given by
\[
X_1=\der_{x_1}-x_2\,\der_{x_3},\quad
X_2=\der_{x_2}+x_1\,\der_{x_3}\quad\mbox{and}\quad
X_3=\der_{x_3}\,.
\]

As mentioned in the introduction, the simplest case of smooth mapping 
that is not Lipschitz is the smooth curve 
$\Gamma:\R\lra\H^1$ defined as $\Gamma(t)=(0,0,t)$.
In this case, using any explicit homogeneous norm, it is straightforward
to see the failure of the Lipschizt property.

On the other hand, if we consider the smooth mapping 
$I:\H^1\lra\H^1$ defined as
\[
I(x_1,x_2,x_3)=(f_1,f_2,f_3)=(x_3,x_1,x_2),
\]
then showing that any of its restrictions to open sets is not locally Lipschitz 
may require slightly more computations using distances and the use of contact equations 
could prove a conveneint tool.
We have
\[
\left\{\begin{array}{l}
f_1X_1f_2-f_2X_1f_1=x_3+x_1x_2 \\
f_1X_2f_2-f_2X_2f_1=-x_1^2
\end{array}\right.
\quad\mbox{and}\qquad
\left\{\begin{array}{l}
X_1f_3=0 \\
X_2f_3=1
\end{array}\right.
\]
Thus, defining $F_1:\H^1\lra\ch$ as
$
F_1(x_1,x_2,x_3)=x_3\,X_1+x_1\,X_2\in\ch
$
and taking into account both \eqref{hcondF} and
\[
[F_1,X_jF_1]=f_1X_jf_2-f_2X_jf_1,
\]
it follows that the contact equations cannot hold at every point.
As a consequence of Theorem~\ref{CharLip}, the mapping $f$
is not Lipschitz on every open subset of $\H^1$.
Clearly, since $I$ is smooth, all components of the mapping $I$
are continuously P-differentiable.
%
%
%
%%%%%%%%%%%%%%%%%%%%%%%%%%%%%%%%%%%%%%%%%%%%%%%%%%%%%%%%%%%%%%%%%%%%%%%%%%%%%%%%%%%%%%%
%
%
%
%
%
%
%
\subsection{Contact equations and Rumin complex}\label{RuminLinear}
%
%
%
%
%
%
%
%%%%%%%%%%%%%%%%%%%%%%%%%%%%%%%%%%%%%%%%%%%%%%%%%%%%%%%%%%%%%%%%%%%%%%%%%%%%%%%%%%%%%%%
%
%
%
%
In view of our study of Lipschitz mappings in the three dimensional Heisenberg group $\H^1$,
we limit ourselves to recall the Rumin complex on $\H^1$, see \cite{Rum90} for the case
of general contact manifolds. We denote by $\Omega^k(\H^1)$ the module of $k$-forms
on $\H^1$ and also
\begin{eqnarray*}
J^2=\left\{\alpha\in\Omega^2(M)\mid \theta\wedge\alpha=0 \right\},\quad
I^1=\left\{\ph\,\theta\mid \ph\in C^\infty(\H^1) \right\}\,,
\end{eqnarray*}
where $\theta=dt+\big(x_2dx^1-x_1dx^2\big)/2$ is the contact form.
In this coordinates, we fix the left invariant vector fields
\[
X_1=\der_1-\frac{x_2}{2}\der_t,\quad X_2=\der_2+\frac{x_1}{2}\der_t\quad
\mbox{and}\quad X_3=\der_t.
\]
We also set
$
\Omega^1(\H^1)/I^1
=\big\{[\alpha_1dx_1+\alpha_2dx_2+\alpha_3\theta]_\cR\mid \alpha_j\in C^\infty(\H^1)\big\}\,.
$
Clearly, 
\[
[\alpha_1dx_1+\alpha_2dx_2+\alpha_3\theta]_\cR=[\alpha_1dx_1+\alpha_2dx_2]_\cR
\]
and we have the following
\begin{The}[Rumin, 1990]
There exists $D:\Omega^1(\H^1)/I^1\lra J^2$ such that
\begin{equation}
0\lra\R\lra C^\infty(\H^1)\lra\Omega^1(\H^1)/I^1\stackrel{D}{\lra} J^2\lra0
\end{equation}
defines a complex whose cohomology coincides with the De Rham cohomology, where
\[
D[\alpha_1dx_1+\alpha_2dx_2]_\cR=d\big(\alpha_1dx_1+\alpha_2dx_2+\alpha_3\theta\big)\in J^2
\quad\mbox{defining}\quad \alpha_3=X_1\alpha_2-X_2\alpha_1.
\]
\end{The}
We denote by $d_\cR$ the differential of this complex.
For more information, we address the reader to \cite{Rum90}.
In the next proposition, we will also use the following notation
\[
A=\left(\begin{array}{cc} a_{11} & a_{12} \\ a_{21} & a_{22} \end{array}\right),
\quad b=\left(\begin{array}{c} b_1 \\ b_2\end{array}\right),\quad
A_j=\left(\begin{array}{c} a_{1j} \\ a_{2j} \end{array}\right)
\quad
\mbox{ and }\quad a=\left(\begin{array}{c} a_1 \\ a_2\end{array}\right).
\]
\begin{Pro}
Let $f=(F_1,f_3):U\subset\H^1\lra\H^1$ be a smooth mapping,
where $U$ is an open neighbourhood of the origin and
$F_1(x)=Ax+at+b$ and $f_3=f_3(x,t)$. Then $f$ is locally Lipschitz if and only if
for some $\tau\in\R$ the following conditions hold
\begin{eqnarray}
&&\det\big(A_1\;a\big)=\det\big(A_2\;a\big)=0\,,\label{Aj}\\
&&f_3(x,t)=\tau+\frac{x_2\det\big(b\;A_2\big)+x_1\det\big(b\;A_1\big)}{2}
+t\left(\frac{\det\big(b\;a\big)}{2}+\det A\right) \label{f_3}\,.
\end{eqnarray}
\end{Pro}
{\sc Proof.}
Since $f$ is Lipschitz and smooth, we have
\begin{eqnarray}\label{CEsyst}
\left\{\begin{array}{l}
X_1f_3=\alpha_1=\frac{1}{2}\left(f_1 X_1f_2-f_2X_1f_1\right) \\
X_2f_3=\alpha_2=\frac{1}{2}\left(f_1 X_2f_2-f_2X_2f_1\right)
\end{array}\right.
\end{eqnarray}
everywhere in $U$. Since $d_\cR$ defines a complex, the Rumin's complex, see \cite{Rum90}, then
\[
d_\cR\left(\big[\alpha_1dx_1+\alpha_2dx_2\big]_\cR\right)=d_\cR\big(d_\cR f_3\big)=0\,.
\]
Taking into account that
\begin{eqnarray*}
d_\cR\big[\alpha_1dx_1+\alpha_2dx_2\big]_\cR&=&
d\big(\alpha_1dx_1+\alpha_2dx_2+\alpha_3\theta\big)\\
&=&\big(X_1\alpha_3-X_3\alpha_1)\,dx_1\wedge\theta+
\big(X_2\alpha_3-X_3\alpha_2)\,dx_2\wedge\theta
\end{eqnarray*}
where $\alpha_3=X_1\alpha_2-X_2\alpha_1$ and $\theta=dt+\big(x_2dx^1-x_1dx^2\big)/2$.
Then the system
\begin{eqnarray}\label{Ruminclosure}
\left\{\begin{array}{l}
X_1\alpha_3-\der_t\alpha_1=0 \\
X_2\alpha_3-\der_t\alpha_2=0
\end{array}\right.
\end{eqnarray}
must hold.
One can check that
\begin{eqnarray*}
&&\!\!\!X_1f_2=a_{21}-\frac{a_2x_2}{2} \quad\qquad  X_1f_1=a_{11}-\frac{a_1x_2}{2} \\
&&\!\!\!X_2f_1=a_{12}+\frac{a_1x_1}{2} \quad\qquad  X_2f_2=a_{22}+\frac{a_2x_1}{2}\,.
\end{eqnarray*}
Thus, a direct calculation shows that
\begin{eqnarray*}
2\alpha_1&=&f_1X_1f_2-f_2X_1f_1 \\
&=&\left(\frac{\det\big(a\;b\big)}{2}-\det(A)\right)x_2+
\frac{\det\big(a\;A_1\big)}{2}x_1x_2+\frac{\det\big(a\;A_2\big)}{2}x_2^2\\
&+&\det\big(a\;A_1\big)t+\det\big(b\;A_1\big)
\end{eqnarray*}
and analogously
\begin{eqnarray*}
2\alpha_2&=&f_1X_2f_2-f_2X_2f_1 \\
&=&\left(\frac{\det\big(b\;a\big)}{2}+\det(A)\right)x_1+
\frac{\det\big(A_2\;a\big)}{2}x_1x_2+\frac{\det\big(A_1\;a\big)}{2}x_1^2\\
&+&\det\big(a\;A_2\big)t+\det\big(b\;A_2\big)\,.
\end{eqnarray*}
Taking into account that $2\alpha_3=X_1(2\alpha_2)-X_2(2\alpha_1)$ and also
\begin{eqnarray*}
&&\!\!\!X_1(2\alpha_2)=-\frac{\det(a\;b)}{2}+\det(A)+\det(A_1\;a)x_1+\det(A_2\;a)x_2\,, \\
&&\!\!\!X_2(2\alpha_1)=-\det(A)+\frac{\det(a\;b)}{2}+\det(a\;A_2)x_2+\det(a\;A_1)x_1
\end{eqnarray*}
we obtain
\[
\alpha_3=\det(A)+x_1\det\left(A_1\;a\right)
+x_2\det\left(A_2\;a\right)-\frac{\det\big(a\;b\big)}{2}\,.
\]
In view of \eqref{Ruminclosure}, we get
\begin{equation}\label{necess}
\det\big(A_1\;a\big)=\det\big(A_2\;a\big)=0\,.
\end{equation}
Notice that by Rumin complex,
this is a necessary condition in order that the system \eqref{CEsyst} admits solutions.
In order to solve the system \eqref{CEsyst}, taking into account \eqref{necess},
direct computations yield
\[
\left\{\begin{array}{l}
2(f_3\circ c_{p,X_1})'(s)=\left(\frac{\det(a\;b)}{2}-\det A\right)p_2+\det\big(b\;A_1\big)\\
2(f_3\circ c_{p,X_2})'(s)=\left(\frac{\det(b\;a)}{2}-\det A\right)p_1+\det\big(b\;A_2\big) 
\end{array}\right.\,
\]
where $c_{p,X_j}(s)=p\exp(sX_j)$.
Set $p=(0,p_2,p_3)$, hence
\[
2f_3\left(p_1,p_2,p_3-\frac{p_1p_2}{2}\right)=2f_3(0,p_2,p_3)+p_1p_2
\left(\frac{\det(a\;b)}{2}-\det A\right)+p_1\det\big(b\;A_1\big)
\]
and $f(0,p_2,p_3)=f_3(0,0,p_3)+p_2\det\big(b\;A_2\big)$.
It follows that
\[
f_3(x)=g\left(t+\frac{x_1x_2}{2}\right)+
\frac{x_2\det\big(b\;A_2\big)+x_1\det\big(b\;A_1\big)}{2}+
\frac{x_1x_2}{2}\left(\frac{\det\big(a\;b\big)}{2}-\det A\right)\,,
\]
where $g:\R\lra\R$ is defined as $g(t)=f_3(0,0,t)$.
By previous computations, we have
\[
\frac{f_1 X_2f_2-f_2X_2f_1}{2}=\frac{x_1}{2}\left(\frac{\det\big(b\;a\big)}{2}+\det A\right)
+\frac{\det\big(b\;A_2\big)}{2}
\]
and our formula for $f_3$ gives
\[
X_2f_3=x_1\,g'\left(t+\frac{x_1x_2}{2}\right)+\frac{\det\big(b\;A_2\big)}{2}
-\frac{x_1}{2}\left(\frac{\det\big(b\;a\big)}{2}+\det A\right)\,.
\]
Imposing the validity of the second equation of \eqref{CEsyst}, it follows that 
\[
g(t)=\tau+\left(\frac{\det\big(b\;a\big)}{2}+\det A\right)t\,.
\]
If $f_3$ satisfies formula \eqref{f_3}, then $f$ satisfies contact equations,
then it is locally Lipschitz. This concludes the proof. $\Box$
%
%
%
%%%%%%%%%%%%%%%%%%%%%%%%%%%%%%%%%%%%%%%%%%%%%%%%%%%%%%%%%%%%%%%%%%%%%%%%%%%%%%%%
%
%
%
%
%
\subsection{Lipschitz mappings as boundary value problems}\label{charextlip}
%
%
%
%
%
%
%%%%%%%%%%%%%%%%%%%%%%%%%%%%%%%%%%%%%%%%%%%%%%%%%%%%%%%%%%%%%%%%%%%%%%%%%%%%%%%%%
%
%
%
This subsection is devoted to the characterization of 
existence of Lipschitz extensions as solutions of contact
equations with assigned boundary datum.
\begin{The}\label{ChaExtLip}
Let $E$ be a closed set of a stratified group $\G$,
let $\Omega=\G\sm E$ and let $f:E\lra\M$ be a Lipschitz mapping.
Existence of a Lipschitz extension $\tilde f:\G\lra\M$
of $f$ with possible larger Lipschitz constant is equivalent to
the existence of a mapping $g:\overline\Omega\lra\M$,
with $G=\exp^{-1}\circ g$, $G_j=\pi_j\circ G$
and $g=\exp\big(G_1+\cdots+G_\upsilon\big)$, such that
\begin{enumerate}
\item
$g_{|\der\Omega}=f_{|\der\Omega}$,
\item
$G_1:\Omega\lra W_1$ is Lipschitz,
\item
all the distributional derivatives $\nabla_HG\in L^\infty_{loc}(\Omega,\cM)^m$
and the system
\begin{equation}\label{hcondG}
\nabla_HG_j=\sum_{n=2}^\upsilon\frac{(-1)^n}{n!}\;\pi_j\left([G,\nabla_HG]_{n-1}\right)
\end{equation}
$\mu$-a.e. holds in $\Omega$, for each $j=2,\ldots\upsilon$.
\end{enumerate}
Whenever a function $g$ satisfying these three conditions exists,
then the Lipschitz extension of $f$ is given by
$\tilde f=f\,{\bf 1}_E+g\,{\bf 1}_\Omega$ and we have the estimate
\begin{eqnarray*}
\Lip(\tilde f)\leq C\;\left(\|\nabla_HG_1\|_{L^\infty(\Omega)}+\Lip(f)	\right),
\end{eqnarray*}
for some geometric constant $C>0$.
\end{The}
{\sc Proof.}
If $f$ admits a Lipschitz extension $\tilde f$,
in view of Theorem~\ref{CharLip}, taking the restriction
$\tilde f_{|\overline{\Omega}}=g$, it follows that $g$
satisfies conditions (1), (2) and (3).
Conversely, let us assume the existence of a mapping $g:\overline\Omega\lra\M$
satisfying these three conditions.
Let $p\in\G$ and let $X\in V_1$ be arbitrarily fixed, with $\|X\|=1$.
Let us consider the curve
\[
\R\ni t\lra c_{p,X}(t)=p\exp(tX)\in \G.
\]
The open set $c_{p,X}^{-1}(\Omega)$ is the disjoint union
$\cup_{j\in\N} I_j$ of open intervals of $\R$, where
$I_j=]a_j,b_j[$. Let $N$ be a homogeneous normal subgroup such that
$\G=N\rtimes H_X$ and let $a_j<a_j'<b_j'<b_j$.
The element $p$ is written in a unique way as
$\overline n\exp(\overline tX)$, due to Proposition~\ref{homdec}.
Then we fix a relatively compact, connected open neighbourhood
$O\subset N$ of the unit element $e$ of $\G$ such that
\[
O\overline n\exp\big((\overline t+s)X\big)\in\Omega\quad\mbox{for every}\quad s\in[a_j',b_j']\,.
\]
By assumptions (2) and (3) on $g$, Lemma~\ref{OJ} applies, hence
we have in particular
\[
\rho\lls g\big(c_{p,X}(b'_j)\big),g\big(c_{p,X}(a'_j)\big)\rls\leq C'\,
\|XG_1\|_{L^\infty(\Omega)}\; d\big(c_{p,X}(b'_j),c_{p,X}(a'_j)\big)\,,
\]
where we have set $C'=C/\inf_{\|Y\|=1} d(\exp Y)$.
Passing to the limit as $a'_j\ra a_j^+$ and $b'_j\ra b_j^-$, we achieve
\begin{equation}\label{a_j}
\rho\lls g\big(c_{p,X}(b_j)\big),g\big(c_{p,X}(a_j)\big)\rls\leq C'\,
\|XG_1\|_{L^\infty(\Omega)}\; d\big(c_{p,X}(b_j),c_{p,X}(a_j)\big).
\end{equation}
Now, consider $\tilde f=f\,{\bf 1}_E+g\,{\bf 1}_\Omega$
and arbitrary select $t,\tau\in\R$, with $t<\tau$.
We have the following cases.
If $t,\tau\in I_j$ for some $j$, then \eqref{a_j} yields
\begin{equation}
\rho\lls \tilde f\big(c_{p,X}(t)\big),\tilde f\big(c_{p,X}(\tau)\big)\rls\leq C'\,
\|XG_1\|_{L^\infty(\Omega)}\; d\big(c_{p,X}(t),c_{p,X}(\tau)\big)\,.
\end{equation}
If $t\in I_j$ and $\tau\notin c_{p,X}^{-1}(\Omega)$, then
the triangle inequality and the fact that $c_{p,X}$ is a geodesic yield
\begin{eqnarray}
&&\rho\lls \tilde f\big(c_{p,X}(t)\big),\tilde f\big(c_{p,X}(\tau)\big)\rls\\
&&\leq C'\,\|XG_1\|_{L^\infty(\Omega)}\; d\big(c_{p,X}(t),c_{p,X}(b_j)\big)
+\Lip(f)\, d\big(c_{p,X}(b_j),c_{p,X}(\tau)\big) \nonumber\\
&&\leq \big(C'\,\|XG_1\|_{L^\infty(\Omega)}+\Lip(f)\big)\,
d\big(c_{p,X}(t),c_{p,X}(\tau)\big)\,. \nonumber
\end{eqnarray}
The same estimate is obtained in the analogous case 
$\tau\in I_j$ and $t\notin c_{p,X}^{-1}(\Omega)$.
If $t\in I_j$ and $\tau\in I_k$ with $j\neq k$, with analogous argument we get
\begin{eqnarray}
\rho\lls \tilde f\big(c_{p,X}(t)\big),\tilde f\big(c_{p,X}(\tau)\big)\rls
\leq \big(2C'\,\|XG_1\|_{L^\infty(\Omega)}+\Lip(f)\big)\,
d\big(c_{p,X}(t),c_{p,X}(\tau)\big)\,. \nonumber
\end{eqnarray}
The remaining case $t,\tau\notin c_{p,X}^{-1}(\Omega)$ is trivial,
since $c_{p,X}(t),c_{p,X}(\tau)\in E$, where $f$ is Lipschitz.
We have shown that for every $r,r'\in c_{p,X}(\R)$, we have
\begin{eqnarray}\label{lipest}
\rho\big(\tilde f(r),\tilde f(r')\big)
\leq \big(2C'\,\|\nabla_HG_1\|_{L^\infty(\Omega)}+\Lip(f)\big)\,
d(r,r')\,. 
\end{eqnarray}
Finally, we adopt the same argument used in the proof of Theorem~\ref{CharLip},
where we connect two arbitrary points $p,p'\in\G$ by a geodesic
with respect to the length distance $\delta$, then we approximate the
geodesic by a sequence of piecewise horizontal lines, according to
Lemma~\ref{HorLinApp}, and we let the estimate \eqref{lipest} pass to the limit. $\Box$

%%%%%%%%%%%%%%%%%%%%%%%%%%%%%%%%%%%%%%%%%%%%%%%%%%%%%%%%%%%%%%%%%%%%%%%%%%%%%%
%%%%%%%%%%%%%%%%%%%%%%%%%%%%%%%%%%%%%%%%%%%%%%%%%%%%%%%%%%%%%%%%%%%%%%%%%%%%%%
%%
%%
%%
%%
%%
\section{Allcock groups}\label{allcockgroups}
%%
%%
%%
%%
%%
%%%%%%%%%%%%%%%%%%%%%%%%%%%%%%%%%%%%%%%%%%%%%%%%%%%%%%%%%%%%%%%%%%%%%%%%%%%%%%
%%%%%%%%%%%%%%%%%%%%%%%%%%%%%%%%%%%%%%%%%%%%%%%%%%%%%%%%%%%%%%%%%%%%%%%%%%%%%%

In this section, we introduce the class of Allcock groups, along with
examples. The characterizing geometric property of these groups 
is related to the notion of isotropic homotopy.
Throughout this section, we fix a 2-step graded algebra $\cn$ with
first layer $\cv$ and second layer $\cz$.
We select a scalar product on $\cn$ such that $\cv$ and $\cz$ are orthonormal.
\begin{Def}{\rm
Let $a,b:[0,1]\lra\cv$ be Lipschitz loops.
We say that $\Gamma:[0,1]^2\lra\cv$ is an {\em isotropic homotopy}
carrying $a$ to $b$ if $(\tau,t)\lra\Gamma(\tau,t)$ is Lipschitz, 
\[
[\der_{\tau}\Gamma,\der_t\Gamma]=0\quad\mbox{a.e. in $[0,1]^2$,}
\]
$\Gamma(\cdot,0)=a$, $\Gamma(\cdot,1)=b$ and $\Gamma(0,\cdot)=\Gamma(1,\cdot)$.}
\end{Def}
\begin{Rem}{\rm
For our purposes, the points $\Gamma(0,\cdot)$ and $\Gamma(1,\cdot)$
need not coincide with some fixed point.}
\end{Rem}
In the sequel, we will use the following linear space
\[
\mbox{Av}_0^\infty=\left\{\sigma\in L^\infty\big(]0,1[,\cz\big)\;\Big|\;
\int_0^1\sigma=0\right\}.
\]
In the next definition the same symbol $|\cdot|$ will
denote both the norm of vectors and 2-vectors. The individual cases
will be clear from the context.
%
%
%
%%%%%%%%%%%%%%%%%%%%%%%%%%%%%%%%%%%%%%%%%%%%%%%%%%%%%%%%%%%%%%%%%%%%%%%%%%%%
%                  ALGEBRAS SURJECTIVE ON ISOTROPIC LOOPS 
%%%%%%%%%%%%%%%%%%%%%%%%%%%%%%%%%%%%%%%%%%%%%%%%%%%%%%%%%%%%%%%%%%%%%%%%%%%%
%
%
%
\begin{Def}\label{suriso}{\rm
We say that $\cn$ is {\em surjective on isotropic loops} if
there exists a constant $C>0$ such that for every $\lambda>0$ 
and every $\sigma\in\mbox{Av}_0^\infty$ one can find a loop
$a\in\Lip([0,1],\cv)$ and an isotropic homotopy $\Gamma:[0,1]^2\lra\cv$ that
carries $a$ to a fixed point, such that $[a,\dot a]=\sigma$ a.e.,
$|a(0)|\leq C\,\lambda$ and the estimates 
\begin{equation}\label{surlipC}
\Lip(\Gamma)\leq C\,\left(\lambda+\frac{\|\sigma\|_{L^\infty}}{\lambda}\right)
\quad\mbox{and}\quad
|\dot a|\leq\frac{C}{\lambda}\,|\sigma|\quad\mbox{a.e.}
\end{equation}
along with 
\begin{equation}\label{surareaC}
\int_0^1\int_0^1|\Gamma_\tau\wedge\Gamma_t|\,d\tau\,dt\leq C\,
\left(\int_0^1|\dot a(t)|\,dt\right)^2\,.
\end{equation}
}\end{Def}
\begin{Rem}{\rm
The $L^\infty$-norm understood for $\sigma\in L^\infty\big(]0,1[,\cz\big)$ 
in the previous definition is given by
\[
\|\sigma\|_{L^\infty}=\max_{j=1,\ldots,s}\|\sigma_j\|_{L^\infty}\,,
\]
where $\sigma=\sum_{j=1}^s\sigma_j\,Z_j$ and $(Z_1,\ldots,Z_s)$
is an orthonormal basis of $\cz$.
}\end{Rem}
%
%
%
%
%%%%%%%%%%%%%%%%%%%%%%%%%%%%%%%%%%%%%%%%%%%%%%%%%%%%%%%%%%%%%%%%%%%%%%%%%
%                        ALLCOCK GROUPS
%%%%%%%%%%%%%%%%%%%%%%%%%%%%%%%%%%%%%%%%%%%%%%%%%%%%%%%%%%%%%%%%%%%%%%%%%%
%
%
%
%
\begin{Def}[Allcock group]
{\rm Let $\cn$ be surjective on isotropic loops
and define
\[
V_1=\bigoplus_{j=1}^n \cv_j,\quad V_2=\cz,\quad [\cv_i,\cv_j]=\{0\}\quad
\mbox{whenever}\quad i\neq j\,,
\]
where all the two step algebras $\cv_j\oplus\zeta$ are
isomorphic to $\cn=\cv\oplus\cz$ for every $j=1,\ldots,n$.
Then the two step algebra $V_1\oplus V_2$ is denoted by $\A\cl_\cn^n$.
This algebra defines a unique stratified group $\A\cl_\cn^n$, that we call
{\em Allcock group of model $\cn$}. If the model $\cn$ is understood,
then we denote an Allcock group simply by $\A\cl^n$.}
\end{Def}
The following example shows that the three dimensional Heisenberg algebra
$\ch^1$ is surjective on isotropic loops.
\begin{Exa}{\rm
Let $\ch$ be the Heisenberg algebra, with layers
$\span\{X_1,X_2\}=\cv$ and $\span\{Z\}=\cz$, where $[X_1,X_2]=Z$.
Let $\sigma\in\mbox{Av}_0^\infty$ and $\lambda>0$, where $\sigma=\sigma_1Z$.
Then we define the Lipschitz curve 
\[
a=a_1X_1+a_2X_2,\quad\mbox{where}\quad
a_1\equiv \lambda\quad\mbox{and}\quad a_2(t)=\frac{1}{\lambda}\int_0^t\sigma_1(s)\,ds\,.
\]
Clearly $[a,\dot a]=\sigma$ and the corresponding homotopy is
\[
\Gamma(\tau,t)=\lambda X_1+(1-t)a_2(\tau)X_2
\] 
is clearly isotropic and carries $a$ to the point $\lambda X_1\in\ch$.
We have $|a(0)|=\lambda$ and simple calculations yield
\[
\Lip(\Gamma)\leq \frac{2\,\|\sigma\|_{L^\infty}}{\lambda}\quad\mbox{and}
\quad |\dot a|=\frac{|\sigma|}{\lambda}\,.
\]
The validity of \eqref{surareaC} is trivial, since $|\Gamma_t\wedge\Gamma_\tau|=0$
at any differentiability point of $\Gamma$.
}\end{Exa}
\begin{Rem}\label{AllcockHeisenberg}{\rm 
The previous example shows that Heisenberg groups $\H^n$ are Allcock groups,
since $\H^n=\A\cl_{\ch}^n$, where $\ch$ is the 3-dimensional Heisenberg algebra.
Let us mention that for Heisenberg groups Lipschitz extensions from the 
Euclidean plane could be also treated by a different method, according 
to F\"assler Master's thesis, 2007.
}\end{Rem}
\begin{Rem}{\rm
Arguing as in the previous example, one can find several examples
of 2-step stratified algebras $\cn$ with one codimensional horizontal
distribution that suitably yield Allcock groups.
}\end{Rem}
On the other hand, it is easy to find Allcock groups where the horizontal 
distribution has codimension higher than one, as in the following 
\begin{Exa}\label{cyclics}{\rm
Let us consider the 2-step algebra $\ck_s=\cv\oplus\cz$, where
$\cv=\{X_1,\ldots,X_{s+1}\}$, $\cz=\span\{Z_1,\ldots,Z_s\}$ and the
only nontrivial bracket relations are
\[
[X_1,X_j]=Z_{j-1},\quad\mbox{and}\quad j=2,\ldots,s+1\,.
\]
Let $\sigma\in\mbox{Av}_0^\infty$ and $\lambda>0$, where 
$\sigma=\sum_{k=1}^s\sigma_kZ_k$.
Then we define the Lipschitz curve 
\[
a=\sum_{j=1}^{s+1}a_jX_j,\quad\mbox{where}\quad
a_1\equiv \lambda\quad\mbox{and}\quad a_j(t)=\frac{1}{\lambda}\int_0^t\sigma_{j-1}(s)\,ds
\]
for every $j=2,\ldots,s+1$. 
Clearly $[a,\dot a]=\sigma$, $|a(0)|=\lambda$ and 
$\dot a=\lambda^{-1}\,\sigma$. Moreover, the mapping
\[
\Gamma(\tau,t)=\lambda X_1+(1-t)\sum_{j=2}^{s+1}a_j(\tau)\,X_j
\] 
is an isotropic homotopy carrying $a$ to $\lambda X_1\in\cv$ and
a simple computation yields
\[
\Lip(\Gamma)\leq\frac{\sqrt{2s}}{\lambda}\,\|\sigma\|_{L^\infty}\,.
\]
Finally, to prove \eqref{surareaC} we observe that there exists $\beta_0>0$
such that
\begin{eqnarray*}
&&\int_{[0,1]^2}|\Gamma_t\wedge\Gamma_\tau|\,d\tau dt
\leq \beta_0\,\int_{[0,1]^2}\bigg|\sum_{j=2}^{s+1}a_jX_j\bigg|\,
\bigg|\sum_{j=2}^{s+1}\dot a_jX_j\bigg|\,d\tau dt\\
&&=\beta_0\,\int_{[0,1]^2}\bigg|\int_0^t\sum_{j=2}^{s+1}\dot a_jX_j\bigg|\,
\bigg|\sum_{j=2}^{s+1}\dot a_jX_j\bigg|\,d\tau dt
\leq\beta_0\,\left(\int_0^1|\dot a|\,dt\right)^2\,.
\end{eqnarray*}
We have proved that $\ck_s$ is surjective on isotrpic loops.
}\end{Exa}
\begin{Rem}{\rm In view of the previous example, we have obtained other
Allcock groups, corresponding to $\A\cl^n_{\ck_s}$.
Notice that they have horizontal distribution of codimension $s$,
for every integer $s\geq1$.
Clearly, a 2-step algebra $\cn$ of $s$-dimensional second layer,
having a subalgebra isomorphic to $\ck_s$ is surjective on isotropic loops.
}\end{Rem}
\begin{Exa}{\rm
We define the ``multi-Heisenberg algebra'' $\cM\ch^s$, as the 2-step
stratified algebra, where first and second layers are spanned 
by the bases $(X_1,\ldots,X_{2s})$ and $(Z_1,\ldots,Z_s)$, respectively,
and the only nontrivial brackets are
\[
[X_j,X_{s+j}]=Z_j\quad\mbox{for every $j=1,\ldots,s.$}
\]
Let $\sigma\in\mbox{Av}_0^\infty$, with $\sigma=\sum_{j=1}^s\sigma_j\,Z_j$
and let $\lambda>0$. We define
\[
a=\lambda\sum_{j=1}^sX_j+\sum_{j=1}^sa_jX_{s+j},\quad\mbox{ where }\quad
a_j(\tau)=\frac{1}{\lambda}\int_0^\tau\sigma_j(s)ds\,.
\]
Then the homotopy $\Gamma(\tau,t)=\lambda\sum_{j=1}^sX_j+(1-t)\sum_{j=1}^sa_j(\tau)X_{s+j}$
is clearly isotropic and carries $a$ to the point $\lambda\sum_{j=1}^sX_j$.
Clearly, $|a(0)|=\lambda\,\sqrt{s}$ and $|\dot a|=\lambda^{-1}\,|\sigma|.$
We also have the estimate
\[
\Lip\big(\Gamma\big)\leq\mbox{ess}\sup_
{\!\!\!\!\!\!\!\!\!\!\!\![0,1]^2}
\bigg(\sum_{j=1}^sa_j^2+\dot a_j^2\bigg)^{1/2}\leq
\frac{\sqrt{2s}}{\lambda}\,\|\sigma\|_{L^\infty}
\]
Finally, arguing exactly as in Example~\ref{cyclics}, we get 
\begin{eqnarray*}
\int_{[0,1]^2}|\Gamma_t\wedge\Gamma_\tau|\,d\tau dt
\leq \beta_0\,\int_{[0,1]^2}\bigg|\sum_{j=1}^sa_jX_{s+j}\bigg|\,
\bigg|\sum_{j=1}^s\dot a_jX_j\bigg|\,d\tau dt
\leq\beta_0\,\left(\int_0^1|\dot a|\,dt\right)^2\,.
\end{eqnarray*}
We have shown that $\cM\ch^s$ is surjective on isotropic loops, 
then $\A\cl^n_{\cM\ch^s}$ are Allcock groups for every $n,s\in\N\sm\{0\}$.
}\end{Exa}
\begin{Exa}{\rm
The {\em complexified Heisenberg algebra} is surjective on isotropic loops.
Recall that this algebra $\C\ch=\cv\oplus\cz$ is an H-type algebra, with
$J_Z:\cv\lra\cv$ and $J_Z^2=-|Z|^2\mbox{I}$ for every $Z\in\cz$.
We fix an orthonormal basis $(Z_1,Z_2)$ of $\cz$ and define the
unit vectors $R_0=X_0$, $R_1=J_{Z_1}X_0$, $R_2=J_{Z_2}X_0$
and $R_3=J_{Z_1}J_{Z_2}X_0$, that form an orthonormal basis of $\cv$.
For more information on the complexified Heisenberg algebra, see \cite{ReiRic}.
Let us fix $\lambda>0$ and choose a curve
$\sigma=\sigma_1 Z_1+\sigma_2Z_2\in\mbox{Av}_0^\infty$.
We define $a=\sum_{j=0}^3a_j\, R_j\in\Lip_0\big([0,1],\cv\big)$ as follows
\[
a_0\equiv \lambda,\quad a_3\equiv0,\quad a_1(t)=\frac{1}{\lambda}\int_0^t\sigma_1
\quad\mbox{and}\quad a_2(t)=\frac{1}{\lambda}\int_0^t\sigma_2\,.
\]
This loop satisfies
$[a,\dot a]=\lambda\,\dot a_1 Z_1+\lambda\,\dot a_2 Z_2=\sigma$,
$|a(0)|=\lambda$ and $|\dot a|=\lambda^{-1}\,|\sigma|$.
We also notice that the isotropic homotopy
$\Gamma(\tau,t)=\lambda\,R_0+(1-t)\big(a_1(\tau)R_1+a_2(\tau)R_2\big)$
carries $a$ to $\lambda R_0\in\cv$ and satisfies $\Lip(\Gamma)\leq 2\,\|\sigma\|_{L^\infty}\,\lambda^{-1}$.
As in the previous cases, the validity of \eqref{surareaC} is straightforward.
This proves that the complexified Heisenberg group $\C\H^1$,
corresponding to $\A\cl_{\C\ch}^1$, is an Allcock group and more generally
$\A\cl_{\C\ch}^n=\C\H^n$ are also Allcock groups. Notice that these groups
are all H-type groups.
}\end{Exa}
\begin{Exa}{\rm Let us consider the {\em quaternionic H-type group}
$\cn_\H$, whose center $\cz$ is spanned by the orthonormal basis
$(Z_1,Z_2,Z_3)$ and $\big(X_0,J_{Z_1}X_0,J_{Z_2}X_0,J_{Z_3}X_0\big)$ is
an orthonormal basis of $\cv$. $X_0$ is a fixed unit vector of $\cv$
and we have
\[
J_{Z_1}J_{Z_2}=J_{Z_3},\quad J_{Z_1}J_{Z_3}=-J_{Z_2}\quad\mbox{and}\quad J_{Z_2}J_{Z_3}=J_{Z_1}\,.
\]
We define the direct product algebra 
$\cn^3_\H=\cv\oplus\cz$, where
$\cv=\cv_1\oplus\cv_2\oplus\cv_3$ 
and $(R_{l0},R_{l1},R_{l2},R_{3l})$ is the orthonormal basis of $\cv_l$ 
for $l=1,2,3$. We have defined  
\[
R_{l0}=X_l,\quad R_{l1}=J_{Z_1}X_l,\quad
R_{l2}=J_{Z_2}X_l, \quad R_{l3}=J_{Z_3}X_l
\]
where $X_l$ is a unit vector of $\cv_l$. Furthermoire, 
whenver $l\neq s$ we set
\[
[R_{li},R_{sj}]=0\quad\mbox{for every $i,j=1,2,3,4$.}
\]
Let $\lambda>0$ and let $\sigma=\sum_{j=1}^3\sigma_j Z_j\in\mbox{Av}_0^\infty$.
We define the curve
\[
a=\lambda\big(\sum_{l=1}^3R_{l0}\big)
+\frac{1}{\lambda}\sum_{l=1}^3\Big(\int_0^t\!\sigma_l(s)\,ds\;\Big) R_{ll}
\in\Lip_0\big([0,1],\cv\big)
\]
Then one can easily check that $[a,\dot a]=\sum_{l=1}^3\sigma_l\,Z_l$,
$|a(0)|=\sqrt{3}\,\lambda$ and $|\dot a|=\lambda^{-1}\,|\sigma|$.
Finally, the isotropic homotopy
\[
\Gamma(\tau,t)=\lambda\big(\sum_{l=1}^3R_{l0}\big)+(1-t)
\frac{1}{\lambda}\sum_{l=1}^3\Big(\int_0^\tau\!\sigma_l(s)\,ds\;\Big) R_{ll}
\]
carries $a$ to $\lambda\big(\sum_{l=1}^3R_{l0}\big)\in\cv$ and satisfies
$\Lip(\Gamma)\leq \sqrt{6}\,\|\sigma\|_{L^\infty}\,\lambda^{-1}.$
Thus, we have proved estimates \eqref{surlipC}. 
Estimates \eqref{surareaC} are obtained as in the previous examples.
We have then proved that the quaternionic H-type group $N^3_\H$
is an Allcock group. Notice that higher dimensional Allcock groups
$\A\cl_{\cn^3_\H}^n$ are H-type groups of dimension $12\,n+3$,
where $n$ is a positive integer.
}\end{Exa}
\begin{Rem}{\rm
Similar computations can be adapted to the octonionic H-type group.
In general, the principle to find Allcock groups is to add as many copies
of the horizontal subspace as possible. This yields the suitable ``room''
to construct isotropic homotopies.
}\end{Rem}
On the other hand, it is not difficult to find 2-step groups
that are not Allcock groups. It suffices to show that some 
2-step stratified algebras are not surjective on isotropic loops,
as we show in the next example.
\begin{Exa}{\rm
The free 2-step free Lie algebra $\cg_{5,2}=V_1\oplus V_2$
on five generators is not surjective on isotropic loops.
Let $(X_1,\ldots,X_5)$ be a basis of generators of $V_1$ and
let $Z_{lp}=[X_l,X_p]$ be the vectors defining a basis of $V_2$
where $1\leq l<p\leq5$. The curve
\[
\sigma(t)=\sum_{1\leq l<p\leq5}\sigma_{lp}(t)\, Z_{lp}
\]
is defined by setting 
\begin{eqnarray}\label{sigmalp}
\sigma_{lp}(t)=\left\{\begin{array}{ll} 
0 & \mbox{if $(l,p)=(1,2)$} \\
t-1/2 & \mbox{if $(l,p)=(1,3)$} \\
t-1/2 & \mbox{if $(l,p)=(2,3)$} \\
t-1/2 & \mbox{if $(l,p)=(1,4)$} \\
t-1/2 & \mbox{if $(l,p)=(2,4)$} \\
t-1/2 & \mbox{if $(l,p)=(1,5)$} \\
0 & \mbox{if $(l,p)=(2,5)$} \\
0 & \mbox{otherwise} 
\end{array}\right.\,.
\end{eqnarray}
clearly belongs to $\mbox{Av}_0^\infty$.
Now, suppose by contradiction that there exists a Lipschitz function $a\in\Lip_0\big([0,1],V_1\big)$ such that 
\[
[a,\dot a]=\sigma\quad\mbox{a.e. in $[0,1]$}\,.
\]
Writing $a(t)=\sum_{j=1}^5 a_j(t)\,X_j$, the previous condition yields
\begin{equation}\label{asigmalp}
\det\left(\begin{array}{cc}
a_l & a_p \\ \dot a_l & \dot a_p
\end{array}\right)=\sigma_{lp}\quad\mbox{a.e. in $[0,1]$}\,.
\end{equation}
Then $\sigma_{12}\equiv0$, $\sigma_{13}=\sigma_{23}\neq0$ a.e. imply
$(a_1,\dot a_1)=\lambda\,(a_2,\dot a_2)$ and $\lambda(t)\neq0$ for
a.e. $t\in[0,1]$. Then constraints \eqref{asigmalp} give $\sigma_{14}=\lambda\,\sigma_{24}$
that yields $\lambda\equiv1$. Therefore the condition
$\sigma_{15}=\lambda\,\sigma_{25}$ corresponds to $\sigma_{15}=\sigma_{25}$,
that conflicts with \eqref{sigmalp}.
}\end{Exa}

%%%%%%%%%%%%%%%%%%%%%%%%%%%%%%%%%%%%%%%%%%%%%%%%%%%%%%%%%%%%%%%%%%%%%%%%%%%%%%
%%%%%%%%%%%%%%%%%%%%%%%%%%%%%%%%%%%%%%%%%%%%%%%%%%%%%%%%%%%%%%%%%%%%%%%%%%%%%%
%%
%%
%%
%%
%%
\section{Lipschitz extensions from the plane to Allcock groups}\label{LipExtSect}
%%
%%
%%
%%
%%
%%%%%%%%%%%%%%%%%%%%%%%%%%%%%%%%%%%%%%%%%%%%%%%%%%%%%%%%%%%%%%%%%%%%%%%%%%%%%%
%%%%%%%%%%%%%%%%%%%%%%%%%%%%%%%%%%%%%%%%%%%%%%%%%%%%%%%%%%%%%%%%%%%%%%%%%%%%%%

This section is devoted to the proof of Theorem~\ref{Allcock1Lipcon}.
Let $\A\cl_\cn^n$ be an Allcock group with $\cn=\cv\oplus\cz$ and
denote by $m$ and $s$ the dimensions of $\cv$ and
$\cz$, respectively. An orthonormal basis $(X_1,\ldots,X_m,Z_1,\ldots,Z_s)$ of $\cn$
will be fixed, where $(X_1,\ldots,X_m)$ and $(Z_1,\ldots,Z_s)$ are bases of $\cv$
and $\cz$, respectively.
We will also choose an orthonormal basis $(X_{ij})_{j=1,\ldots,m}$
of $\cv_i$ and an orthonormal basis $(Z_j)_{j=1,\ldots,s}$ of $\cz$.
Since $\cv_j$'s are all isomorphic to $\cv$, then we can select $X_{ij}$
such that 
\begin{equation}\label{isomprod}
[X_{il},X_{ip}]=[X_l,X_p]\quad\mbox{for every}\quad 
i=1,\ldots,n\quad\mbox{and}\quad 1\leq l<p\leq m\,.
\end{equation}
Then we fix graded coordinates in $\A\cl_\cn^n$ with respect to this basis.
We consider
\[
(x,y)=\sum_{\substack{i=1,\ldots,n\\j=1,\ldots,m}}
x_{ij}\,e_{ij}+\sum_{l=1,\ldots,s}y_l\,E_l
\]
where $(e_{ij},E_l)$ is the canonical basis of $\R^{mn}\times\R^s$.
Precisely, a point in $\A\cl_\cn^n$ of coordinates $(x,y)$ is given by
\[
\exp\Bigg(\sum_{\substack{i=1,\ldots,n\\j=1,\ldots,m}}
x_{ij}\,X_{ij}+\sum_{l=1,\ldots,s}y_l\,Z_l\Bigg)\in\A\cl_\cn^n\,,
\]
where $\exp:\cA\cl_\cn^n\lra\A\cl_\cn^n$ is the standard exponential mapping.
Due to \eqref{isomprod}, for every $i=1,\ldots,n$, we have
\[
[X_{il},X_{ip}]=[X_l,X_p]=\sum_{k=1}^s b_{lp}^k\; Z_k\,.
\]
Then, we are in the position to introduce the {\em multi-symplectic form}
\begin{equation}\label{multsimpl}
\omega=\sum_{k=1}^s\omega^k\; E_k,\quad\mbox{where}\quad
\omega^k=\sum_{j=1}^n\,\sum_{1\leq l<p\leq m}\;b_{lp}^k\;\;
dx_{jl}\wedge dx_{jp}\,.
\end{equation}
Notice that $\omega$ is an $\R^s$-valued 2-form defined on $\R^{mn}$.
A primitive of this form is 
\begin{equation}\label{multprim}
\theta=\sum_{k=1}^s\theta^k\; E_k,\quad\mbox{where}\quad
\theta^k=\sum_{j=1}^n\,\sum_{1\leq l<p\leq m}\frac{b_{lp}^k}{2}\;\;
\big(x_{jl}\,dx_{jp}-x_{jp}\,dx_{jl}\big).
\end{equation}
Then for every Lipschitz curve $c:[a,b]\lra\R^{mn}$, we define
the {\em multi-symplectic area} swept by $c$ as follows
\[
\int_c\theta=\sum_{k=1}^s\bigg(\int_a^b\,(c^*\theta^k)(t)\;dt\bigg)\,E_k\in\R^k\,.
\]
This definition can be clearly extended to one-dimensional compact Lipschitz
manifolds through local parametrizations.
\begin{Def}{\rm
Let $\A\cl_\cn^n$ be an Allcock group and let $\R^{mn}\times\R^s$
be the associated graded coordinates, where $\R^{mn}$ is the space
of coordinates of the first layer. Then $\omega$ and $\theta$ defined
above are the {\em associated} multi-simplectic form and its primitive,
respectively.}
\end{Def}
\begin{Def}{\rm
If $c_i:[0,1]\lra\R^{mn}$, $i=0,1$, are Lipschitz loops,
we say that $\Gamma:[0,1]^2\lra\R^{mn}$ is a {\em multi-isotropic homotopy}
carrying $c_0$ to $c_1$ if $\Gamma$ is Lipschitz, the pull-back
$\R^s$-valued 2-form $\Gamma^*\omega$ a.e. vanishes in $[0,1]^2$,
$\Gamma(\cdot,0)=c_0$, $\Gamma(\cdot,1)=c_1$
and $\Gamma(0,\cdot)=\Gamma(1,\cdot)$.}
\end{Def}
All the preceeding notions will play a key role in the proof of
the next theorem. 
Recall that the Euclidean norm will be understood on $\R^{mn}$ and $S^1$ will be thought
of as the subset $\{(x,y)\mid x^2+y^2=1\}$ equipped with the Euclidean 
distance of $\R^2$. The symbol $D$ denotes the closed unit disk of $\R^2$.
\begin{The}\label{allcock}
Let $\A\cl^n$ be an Allcock group and let $\omega$ and $\theta$
be the associated forms on $\R^{mn}$ with respect to graded coordinates,
wehre $n\geq2$. Then there exists a geometric constant $\kappa>0$ such
that for every Lipschitz loop $c:S^1\lra\R^{mn}$ with $\int_c\theta=0$
and $c(1,0)=0$, then there exists a Lipschitz extension $\ph:D\lra\R^{mn}$
such that $\ph^*(\omega)=0$ a.e. in $D$ and there exists 
\begin{equation}
\Lip(\ph)\leq \kappa\;\Lip(c)\,.
\end{equation}
\end{The}
{\sc Proof.}
We consider $\alpha(t)=c\big(e^{2\pi i t}\big)$ and 
set $\tilde\alpha:\R\lra\R^{mn}$ such that $\tilde\alpha_{|[0,1]}=\alpha$ 
and $\tilde\alpha(\tau)=\alpha(0)=\alpha(1)=0$ for every $\tau\in\R\sm[0,1]$.
Clearly we have
\[
\Lip(\alpha)=\Lip(\tilde\alpha)\leq 2\pi\,\Lip(c)\,.
\]
We first change the parametrization by the homotopy
\[
\Gamma_1:[0,1]^2\lra\R^{mn},\;
\Gamma_1(\tau,t)=\tilde\alpha\big((1+t)\tau-t\big)\,.
\]
Notice that $\Gamma_1$ is Lipschitz and clearly the pull-back form
$\Gamma_1^*\omega$ vanishes a.e. in $[0,1]^2$.
We have shown that $\alpha$ is isotropically homotopic to
$[0,1]\ni \tau\ra\tilde\alpha(2\tau-1)$.
Now, we write $\tilde\alpha(\tau)=\beta_1(\tau)+\beta_2(\tau)$,
where we have set
\begin{eqnarray*}
\beta_1(\tau)=\sum_{j=1}^m \tilde\alpha_{1j}(\tau)\,e_{1j}\quad\mbox{and}\quad
\beta_2(\tau)=\sum_{\substack{i=2,\ldots,n\\j=1,\ldots,m}}
\tilde\alpha_{ij}(\tau)\,e_{ij}\,.
\end{eqnarray*}
Taking into account formula \eqref{multsimpl} for the
multi-simplectic 2-form $\omega$,
we notice that
\[
\omega\big(\beta_1(\tau)\wedge\beta_2(\tau)\big)=
\sum_{k=1}^s\omega^k\big(\beta_1(\tau)\wedge\beta_2(\tau)\big)\,E_k=0,
\]
since $dx_{jl}\wedge dx_{jp}\big(e_{1r}\wedge e_{iu}\big)=0$
whenever $i\neq 1$.
Therefore, defining the homotopy
\[
\Gamma_2:[0,1]^2\lra\R^{mn},\;
\Gamma_2(\tau,t)=\beta_1(2\tau-1+t)+\beta_2(2\tau-1),
\]
it follows that for a.e. $(\tau,t)\in[0,1]^2$, we have
\begin{eqnarray*}
\Gamma_2^*\omega=\omega\big(\der_{\tau}\Gamma_2\wedge\der_t\Gamma_2\big)=
2 \omega\big(\dot\beta_2(2\tau-1)\wedge\dot\beta_1(2\tau-1+t)\big)=0\,.
\end{eqnarray*}
Then $\tilde\alpha(2\tau-1)$ is multi-isotropically homotopic to the product curve
\[
(\beta_1\beta_2)(\tau)=\beta_1(2\tau)+\beta_2(2\tau-1)\,.
\]
Now, we move $\beta_1$ to a multi-symplectically orthogonal space defining
\[
\overline{\beta}_1(\tau)=\sum_{j=1}^m\tilde\alpha_{1j}(\tau)\,e_{2j}\,.
\]
Thus, we get $\omega(\beta_1\wedge\overline\beta_1)=0$. Notice that to perform
this move from $\beta_1$ to $\overline\beta_1$, we have used the assumption $n\geq2$.
We define the isotropic homotopy $\Gamma_3:[0,1]^2\lra\R^{mn},$
\[
\Gamma_3(\tau,t)=\cos(t\pi/2)\,\beta_1(2\tau)+\sin(t\pi/2)\,\overline\beta_1(2\tau)+
\beta_2(2\tau-1)\,.
\]
Now, we observe that 
\begin{eqnarray*}
\omega^k(\beta_1\wedge\dot\beta_1)&=&\sum_{j=1}^n\,
\sum_{1\leq l<p\leq m}\;b_{lp}^k\;\;(dx_{jl}\wedge dx_{jp})(\beta_1\wedge\dot\beta_1)\\
&=&\sum_{1\leq l<p\leq m}\;b_{lp}^k\;\;(dx_{1l}\wedge dx_{1p})(\beta_1\wedge\dot\beta_1)
\end{eqnarray*}
and similarly we have
\begin{eqnarray*}
\omega^k(\overline\beta_1\wedge\dot{\overline\beta}_1)&=&\sum_{j=1}^n\,
\sum_{1\leq l<p\leq m}\;b_{lp}^k\;\;(dx_{jl}\wedge dx_{jp})
(\beta_1\wedge\dot{\overline\beta}_1)\\
&=&\sum_{1\leq l<p\leq m}\;b_{lp}^k\;\;(dx_{2l}\wedge dx_{2p})
(\overline\beta_1\wedge\dot{\overline\beta}_1)\,,
\end{eqnarray*}
hence $\omega(\beta_1\wedge\dot\beta_1)
=\omega(\overline\beta_1\wedge\dot{\overline\beta}_1\big)$ a.e. in $[0,1]$.
Thus, a simple computation yields
\begin{eqnarray*}
\Gamma_3^*\omega&=&\omega\big(\der_t\Gamma_3\wedge\der_\tau\Gamma_3\big)\\
&=&-\pi\sin(t\pi/2)\cos(t\pi/2)\,
\Big(\omega\big(\beta_1(2\tau)\wedge\dot\beta_1(2\tau)\big)-
\omega\big(\overline\beta_1(2\tau)\wedge\dot{\overline\beta}_1(2\tau)\big)\Big)\\
&&+\cos(t\pi/2)\,\omega\big(\overline\beta_1(2\tau)\wedge\dot\beta_2(2\tau-1)\big)\\
&=&\cos(t\pi/2)\,\omega\big(\overline\beta_1(2\tau)\wedge\dot\beta_2(2\tau-1)\big)
\end{eqnarray*}
for every $t\in[0,1]$ and a.e. $\tau\in[0,1]$.
On the other hand, $\overline\beta_1(2\tau)$ vanishes when 
$2^{-1}\leq \tau\leq1$ and $\dot\beta_2(2\tau-1)=0$ for a.e. $\tau\in[0,2^{-1}]$, this
implies that $\Gamma_3^*\omega=0$ a.e. in $[0,1]^2$.
We have proved that $\beta_1\beta_2$ is multi-isotropically homotopic to
\[
(\overline\beta_1\beta_2)(\tau)=\overline\beta_1(2\tau)+\beta_2(2\tau-1)\,.
\] 
To construct the next homotopy, we use the fact that the algebra
$\cn$ is surjective on isotropic loops.
Let us consider the $L^\infty$ curve
\begin{eqnarray}\label{sigmaomega}
\sigma=\sum_{k=1}^s\,\omega^k\Big(\big(\overline\beta_1\beta_2\big)\wedge
\frac{d}{ds}\big(\overline\beta_1\beta_2\big)\Big)\,Z_k\,.
\end{eqnarray}
We wish to show that $\sigma\in\mbox{Av}_0^\infty$, namely
$\int_0^1\sigma=0$. To do so, we use the vanishing of the multi-simplectic
area $\int_{c}\theta=0$, that yields
\[
\int_{S^1}c^*\theta=\int_0^1\alpha^*\theta=0\,.
\]
Then we get 
\[
\int_0^1 (\overline\beta_1\beta_2)^*\theta
=\int_0^{1/2}(\overline\beta_1)^*\theta
+\int_{1/2}^1\beta_2^*\theta=\int_0^{1/2}\beta_1^*\theta
+\int_{1/2}^1\beta_2^*\theta=
\int_0^1(\beta_1\beta_2)^*\theta=\int_0^1\alpha^*\theta=0.
\]
and the equalities 
\[
(\overline\beta_1\beta_2)^*\theta=\frac{1}{2}\;
\omega\Big(\big(\overline\beta_1\beta_2\big)\wedge
\frac{d}{ds}\big(\overline\beta_1\beta_2\big)\Big)=\frac{1}{2}\;\sigma
\]
prove our claim, namely, $\sigma\in\mbox{Av}_0^\infty$.
Since $\cn$ is surjective on isotropic loops, we set
\begin{equation}\label{lambdaL}
\lambda=\mbox{length}(c)=\int_0^1|\dot\alpha(t)|\,dt\leq 2\pi\,\Lip(c)=2\pi L
\end{equation}
and apply Definition~\ref{suriso}.
Then we can find a loop $a\in\Lip([0,1],\cv)$ 
that is isotropically homotopic to a point and in particular satisfies
\begin{equation}\label{aLipC}
[a,\dot a]=\sigma\qquad\mbox{and}\qquad
\max_{[0,1]}|a|\leq C\left(\lambda+\frac{\|\sigma\|_{L^\infty}}{\lambda}
\right)\,.
\end{equation}
We write $a=\sum_{i=1}^ma_i\,X_i$, then
\[
\sum_{1\leq l<p\leq m}\,\det\left(\begin{array}{cc}
a_l & a_p \\ \dot a_l & \dot a_p
\end{array}\right)\,b_{lp}^k=\omega^k\Big(\big(\overline\beta_1\beta_2\big)\wedge
\frac{d}{ds}\big(\overline\beta_1\beta_2\big)\Big)
\]
for every $k=1,\ldots,s$. We introduce the Lipschitz loop
$q=\sum_{j=1}^ma_j\,e_{1j}$ in $\R^{mn}$ and observe that
\[
\omega\big(q(\tau)\wedge(\overline\beta_1\beta_2)(\tau)\big)=0\,.
\]
Furthermore, by definition of $\omega^k$, one immediately gets 
\begin{equation}\label{alap}
\omega^k\big(q\wedge\dot q\big)=\sum_{1\leq l<p\leq m}\,\det\left(\begin{array}{cc}
a_l & a_p \\ \dot a_l & \dot a_p
\end{array}\right)\,b_{lp}^k\,.
\end{equation}
It follows that
\[
\omega\big(q\wedge\dot q\big)=
\omega\big((\overline\beta_1\beta_2)\wedge
\frac{d}{d\tau}(\overline\beta_1\beta_2)\big)
\quad\mbox{a.e. in $[0,1]$}\,.
\]
Due to the previous condition, arguing as for $\Gamma_3$, one 
easily finds that the Lipschitz mapping $\Gamma_4:[0,1]^2\lra\R^{mn}$ defined as 
\[
\Gamma_4(\tau,t)=\cos(t\pi/2)\,(\overline\beta_1\beta_2)(\tau)
+\sin(t\pi/2)\,q(\tau)
\]
is a multi-isotropic homotopy between $\overline\beta_1\beta_2$ and $q$.
According to Definition~\ref{suriso}, there exists an isotropic homotopy
$\Gamma:[0,1]^2\lra\cv$ such that $\Gamma(\cdot,0)=a$,
$\Gamma(\cdot,1)\equiv\xi$, $\Gamma(0,\cdot)=\Gamma(1,\cdot)$ and
$[\der_{\tau}\Gamma,\der_t\Gamma]=0$ a.e. in $[0,1]^2$.
Writing $\Gamma=\sum_{i=1}^n\sum_{j=1}^m\Gamma_{ij}X_{ij}$ and defining 
$\Gamma_5:[0,1]^2\lra\R^{mn}$ as $\Gamma_5=\sum_{i=1}^n\sum_{j=1}^m\Gamma_{ij}e_{ij}$, 
we have that $\Gamma_5^*\omega=0$ and estimate \eqref{surlipC} yields
a constant $\tilde C_1>0$ such that
\begin{equation}\label{LipGamma5A}
\Lip(\Gamma_5)\leq \tilde C_1\, \left(\lambda+\frac{\|\sigma\|_{L^\infty}}
{\lambda}\right)\,.
\end{equation}
Defining $\xi=\sum_{i,j}\xi_{ij}X_{ij}$ and $q_0=\sum_{i,j}\xi_{ij}e_{ij}\in\R^{mn}$,
one immediately check that $\Gamma_5$ makes $q$ multi-isotropically
homotopic the a point $q_0\in\R^{mn}$. 
Thus, pasting all the previous multi-isotropic homotopies, we get the
following mapping $H:[0,1]^2\lra\R^{mn}$, defined as
\[
H(\tau,t)=\left\{\begin{array}{lr}
\Gamma_1(\tau,5t)   &  0\leq t<1/5 \\ 
\Gamma_2(\tau,5t-1) &  1/5\leq t<2/5 \\ 
\Gamma_3(\tau,5t-2) &  2/5\leq t<3/5 \\ 
\Gamma_4(\tau,5t-3) &  3/5\leq t<4/5 \\
\Gamma_5(\tau,5t-4) &  4/5\leq t<1 
\end{array}\right.\,.
\]
It is clearly both continuous and a.e. differentiable in $[0,1]^2$.
Moreover, we also have $H^*(\omega)=0$ a.e.in $[0,1]^2$.
By convexity of $[0,1]^2$ and triangle inequality, one easily notices 
the following estimate
\begin{equation}\label{max5}
\Lip(H)\leq 5\;\sum_{j=1}^5\Lip(\Gamma_j)\,.
\end{equation}
Direct computations show that
\begin{eqnarray}\label{LipComp}
\left\{\begin{array}{l}
\Lip(\Gamma_1)\leq 6\pi\,L\\
\Lip(\Gamma_2)\leq 6\pi\,L \\
\Lip(\Gamma_3)\leq 6\pi L+\frac{\pi}{2}\max_{[0,1]}|\alpha|
\end{array}\right.\,,
\end{eqnarray} 
where we recall that we have defined $L=\Lip(c)$.
Taking into account the definition of Allcock group, we have
to estimate $\|\sigma\|_{L^\infty}$. 
Thus, taking into account \eqref{sigmaomega}, a direct computation yields 
\[
|\omega^k(\gamma(\tau)\wedge\dot\gamma(\tau))|\leq 2\,n\,m\,
\bigg(\!\max_{\substack{r=1,\ldots,s\\ 1\leq l<p\leq m}}|b_{lp}^r|\,\bigg)
|\gamma(\tau)|\,|\dot\gamma(\tau)|
\]
for every $k=1,\ldots,s$, hence estimates
\[
\max_{[0,1]}\left|\big(\overline{\beta_1}\beta_2\big)\right|
\leq\max_{[0,1]}|\alpha|\quad
\mbox{and}\quad
\left\|\frac{d}{d\tau}\big(\overline{\beta_1}\beta_2\big)\right\|
_{L^\infty\big((0,1),\R^{mn}\big)}\leq2\|\dot\alpha\|_{L^\infty\big((0,1),\R^{mn}\big)}
\leq 4\pi L\,,
\]
along with $\max_{[0,1]}|\alpha|\leq\,\lambda$, lead us to the estimate
\[
\|\sigma\|_{L^\infty}\leq 8\pi
\sqrt{s}\,n\,m\,\bigg(\!\max_{\substack{r=1,\ldots,s\\ 1\leq l<p\leq m}}
|b_{lp}^r|\,\bigg)\,\lambda\,L=C_0\lambda L\,.
\]
As a consequence, by virtue of \eqref{lambdaL} and \eqref{LipGamma5A},
we get a constant $C_1>0$ such that
\begin{equation}\label{LipGamma5}
\Lip(\Gamma_5)\leq C_1\,L\,.
\end{equation}
In addition, by \eqref{lambdaL} and \eqref{aLipC}, observing that $\max|a|=\max|q|$, 
we also obtain
\[
\max_{[0,1]}|q|\leq C\left(2\pi+C_0\right)\,L\,,
\]
where $C_0>0$ is a geometric constant depending on the group.
It is also obvious that $\Lip(q)\leq\Lip(\Gamma_5)\leq C_1\,L$.
Thus, joining all the previous estimates, we get a new constant $C_2>0$
such that 
\[
\Lip(\Gamma_4)\leq C_2\;L\,.
\]
By last inequality along with \eqref{LipGamma5},
all previous estimates for $\Lip(\Gamma_j)$ and applying \eqref{max5}, 
we have found a geometric constant $\kappa_1>0$,
only depending on the group, such that
\begin{equation}\label{LipH}
\Lip(H)\leq \kappa_1\;L\,.
\end{equation}
Then $H$ is a multi-isotropic homotopy such that
$H(\cdot,0)=\alpha:[0,1]\lra\R^{mn}$ and $H(\cdot,1)\equiv q_0$.
The condition $H(0,t)=H(1,t)$ for every $t\in[0,1]$ implies that
\[
\ph\left(\rho\,e^{i2\pi\theta}\right)=\left\{\begin{array}{ll}
H\big(\theta,2(1-\rho)\big)  &
 \mbox{if}\; 1/2\leq \rho\leq 1 \\
q_0  & \mbox{if}\; 0\leq \rho\leq1/2 \end{array}\right.
\]
is well defined on the closed unit disk $D\subset\R^2$.
Furthermore, a direct computation shows that 
$\Lip(\ph_{D\sm B_{1/2}})\leq c_0\,\Lip(H)$,
for a suitable geometric constant $c_0>0$, where 
$B_{1/2}=\{z\in\R^2\mid |z|<1/2\}$.
Since $\ph_{|\overline{B}_{1/2}}\equiv q_0$, then
\[
\Lip(\ph)\leq c_0\,\Lip(H)\leq c_0\,\kappa_1\,L
\]
Clearly $\ph_{S^1}\equiv c$ and $\ph^*\omega=0$ a.e. in $D$.
Our claim is achieved. $\Box$
\begin{Rem}{\rm
The previous theorem extends Theorem~2.3 of \cite{All}, where 
the standard simplectic space $\R^{2n}$ is replaced with the 
$\R^{mn}$ equipped with the multi-simplectic form $\omega=\sum_{k=1}^s\,\omega^k\,E_k$.
}\end{Rem}
Next, we will show how Theorem~\ref{allcock} leads us to a Lipschitz extension theorem.
We will use some abstract tools in metric spaces, following
the work by Lang and Schlichenmaier, \cite{LanSch}.
\begin{Def}{\rm
We say that a metric space $Y$ is {\em Lipschitz $m$-connected}
for some $m\in\N$ if there exists a constant $c_m>0$ such that 
any Lipschitz map $\Gamma:S^m\lra Y$ has 
a Lipschitz extension $\Phi:D^{m+1}\lra Y$ with estimate 
$\Lip(\Phi)\leq c_m\,\Lip(\Gamma
)\,.$
}\end{Def}
\begin{Def}{\rm
Let $(X,Y)$ be a couple of metric spaces.  We say that $(X,Y)$ has the
{\em Lipschitz extension property} if there exists $C>0$ such that
for every subset $Z\subset X$ and every Lipschitz map $f:Z\lra Y$,
there exists a Lipschitz extension $\overline f:X\lra Y$ such that
$\Lip(\overline f)\leq C\,\Lip(f).$
}\end{Def}
\begin{The}[Lang-Schlichenmaier, \cite{LanSch}]\label{LanSch}
Let $X$ and $Y$ be two metric spaces and suppose that the Nagata
dimension of $X$ is less than or equal to $n$ and that $Y$ is complete.
If $Y$ is Lipschitz $m$-connected for $m=0,1,\ldots,n-1$, then
the pair $(X,Y)$ has the Lipschitz extension property.
\end{The}
Taking into account that Nagata dimension of $\R^2$ is clearly two,
by Theorem~\ref{LanSch} it follows that both Lipschitz 0-connectedness and 
Lipschitz 1-connectedness of $\A\cl_\cn^n$ imply Corollary~\ref{LipExtR2Aln}. 
The former property is a consequence of the fact that $\A\cl_\cn^n$, 
as any stratified group, is connected by geodesics.
The latter is proved in the following
\begin{The}\label{1Lipcon}
$\A\cl^n$ is 1-connected for every $n\geq2$.
\end{The}
{\sc Proof.}
Let $\Gamma:S^1\lra\A\cl^n$ be a Lipschitz loop.
Up to a left translation, that preserves the Lipschitz constant of $\Gamma$,
we can assume that 
$\Gamma(1,0)=e$, where $e$ is the unit element of $\A\cl^n$.
We introduce the 1-periodic mappings
$a:[0,1]\lra V_1$ and $b:[0,1]\lra V_2$ such that
\[
\Gamma\big(e^{2\pi it}\big)=\exp\big(a(t)+b(t)\big)\,.
\]
We consider our fixed basis $(X_{ij})$ of $V_1$, that satisfies \eqref{isomprod},
along with the orthonormal basis $(Z_k)$ of $V_2$.
We define
\[
a(t)=\sum_{\substack{1\leq i\leq n\\ 1\leq j\leq m}}
\alpha_{ij}(t)\,X_{ij}\quad\mbox{and}\quad b(t)=\sum_{k=1}^s\beta_k(t)\,Z_k
\]
where $\alpha=(\alpha_{ij}):[0,1]\lra\R^{mn}$ and $\beta=(\beta_k):[0,1]\lra\R^s$
are Lipschitz loops that satisfy $\alpha(0)=\alpha(1)=0$ and $\beta(0)=\beta(1)=0$.
Since $\Gamma$ is Lipschitz, Theorem~\ref{CharLip} implies the a.e. validity
of contact equations
\begin{equation}\label{ce}
\sum_{k=1}^s\dot\beta_k\,Z_k=\frac{1}{2}\;[a,\dot a]
=\frac{1}{2}\sum_{j=1}^n\sum_{1\leq l<p\leq m}
\sum_{k=1}^s b_{lp}^k\,\big(\alpha_{jl}\dot\alpha_{jp}-
\alpha_{jp}\dot\alpha_{jl}\big)\,Z_k=\sum_{k=1}^s\alpha^*(\theta^k)Z_k
\end{equation}
As a consequence, the 1-periodicity of $\beta$ yields
\begin{equation}\label{vansymparea}
\int_\alpha\theta=\sum_{k=1}^s\bigg(\int_0^1\,(\alpha^*\theta^k)(t)\;dt\bigg)\,E_k=0\,.
\end{equation}
We define the curve $c=(c_{ij}):S^1\lra\R^{mn}$ defined by
\[
c_{ij}\big(e^{2\pi i t}\big)=\alpha_{ij}(t)\,,
\]
therefore $c$ has vanishing multi-symplectic area, due to \eqref{vansymparea},
and $c(1,0)=0$.
This allows us to apply Theorem~\ref{allcock}, getting a Lipschitz
extension $\ph:D\lra\R^{mn}$ of $c:S^1\lra\R^{mn}$ such that
\begin{equation}\label{ph*omega0}
\ph^*\omega=0\quad\mbox{a.e. in $D$}\,.
\end{equation}
To construct the extension of $\Gamma$, we introduce the function
\[
A:D\lra V_1,\qquad A(x_1,x_2)=\sum_{\substack{1\leq i\leq n\\ 1\leq j\leq m}}\ph_{ij}(x_1,x_2)\,X_{ij}\,.
\]
We observe that $a(t)=A\big(e^{2\pi i t}\big)$ and define 
$z_k(e^{2\pi i t})=\beta_k(t)$, where $z_k:S^1\lra\R^s$ is Lipschitz continuous
for every $k=1,\ldots,s$.
To achieve our claim, we have to find a Lipschitz extension
$T:D\lra V_2$ such that $T_{|S^1}=\sum_{k=1}^sz_k\,Z_k$ and the following
contact equations a.e. hold
\begin{equation}\label{contZ}
\der_{x_1}T=\frac{1}{2}\big[A,\der_{x_1}A]\quad\mbox{and}\quad
\der_{x_2}T=\frac{1}{2}\big[A,\der_{x_2}A]\,.
\end{equation}
In fact, if such a function $B$ exists,
then applying Theorem~\ref{RMLip} in the case $N$ is closed Euclidean disck $D$ of $\R^2$,
hence the mapping
\[
\Phi=\exp\circ(A+T):D\lra\A\cl^n
\]
is Lipschitz continuous and 
\[
\Lip(\Phi)\leq C\,\Lip(A)=C\,\Lip(\ph),
\]
for a suitable geometric constant.
Furthermore, $\Phi_{|S^1}=\Gamma$, hence Theorem~\ref{allcock} yields another
geometric constant $\kappa$ such that
\begin{equation}\label{LipPhi}
\Lip(\Phi)\leq \kappa\;C\,\Lip(c)\,.
\end{equation}
By Lemma~\ref{LipF_1}, there exists $C_1>0$ such that
\begin{equation}\label{LipPhiGamma}
\Lip(\Phi)\leq \kappa\; C\;C_1\;\Lip(\Gamma)\,.
\end{equation}
Thus, to conclude the proof, we are left to show the existence of
$T:D\lra V_2$ satisfying the contact equations
\eqref{contZ} and the boundary condition $T_{|S^1}=\sum_{k=1}^sz_k\,Z_k$.

First, for every $k=1,\ldots,s$, we consider the 1-form
$f_k=f_k^1dx_1+f_k^2dx_2$ of components
\[
f_k^i=\frac{1}{2}\sum_{j=1}^{n}\sum_{1\leq l<p\leq m} b_{lp}^k\;
\big(\ph_{jl}\;\der_{x_i}\ph_{jp}-\ph_{jp}\;\der_{x_i}\ph_{jl}\big)\in L^\infty(B_1),
\]
where $B^1$ is the unit open ball of $\R^2$.
Smoothing the function $\ph$ with $\ph^\ep=\ph*\zeta_\ep$,
we obtain the approximating forms
$g_{k\ep}=g_{k\ep}^1dx_1+g_{k\ep}^2dx_2$ of components
\[
g_{k\ep}^i=\frac{1}{2}\sum_{j=1}^{n}\sum_{1\leq l<p\leq m} b_{lp}^k\;
\big(\ph_{jl}^\ep\;\der_{x_i}\ph_{jp}^\ep-\ph_{jp}^\ep\;\der_{x_i}\ph_{jl}^\ep\big)\,,
\]
where $\zeta_\ep$ is a standard mollifier and $g_{k\ep}$ is a smooth
1-form on $B_{1-\ep}$. It follows that
\begin{eqnarray*}
dg_{k\ep}&=&\big(\der_{x_1}g_{k\ep}^2-\der_{x_2}g_{k\ep}^1\big)dx_1\wedge dx_2\\
&=&\sum_{j=1}^n\sum_{1\leq l<p\leq m
}b_{lp}^k\;\big(\der_{x_1}\ph^\ep_{jl}\;\der_{x_2}\ph^\ep_{jp}-\der_{x_1}\ph^\ep_{jp}
\der_{x_2}\ph^\ep_{jl}\big)\;dx_1\wedge dx_2\\
&=&(\ph^\ep)^*\omega^k\,.
\end{eqnarray*}
Furthermore, $(\ph^\ep)^*\omega^k$ a.e. converge to $\ph^*\omega^k$ and has 
uniformly bounded $L^\infty$-norm with respect to $\ep$. 
Thus, taking into account \eqref{ph*omega0},
we have proved that $(\ph^\ep)^*\omega^k$ converges to zero in 
$L^1(B_{1-\delta})$ as $\ep\ra0^+$ for arbitrary $0<\delta<1$. 
Taking also into account that $g_{k\ep}$ a.e. converge to $f_k$
and are uniformly bounded with respect to $\ep$ in the $L^\infty$-norm,
it follows that
\[
df_k=0\quad\mbox{in the distributional sense.}
\]
As a consequence, setting $f_{k\ep}^i=f_k^i*\zeta_\ep$ 
and $f_{k\ep}=f_{k\ep}^1 dx_1+f_{k\ep}^2dx_2$ one gets
\[
df_{k\ep}=0\quad\mbox{on $B_{1-\ep}$}\,.
\]
This gives the existence of a unique function $\tilde\psi_{k\ep}\in C^\infty(B_{1-\ep})$
such that $\tilde\psi_{k\ep}(0,0)=0$ and $d\tilde\psi_{k\ep}=f_{k\ep}$ in $B_{1-\ep}$\,.
The family of functions $\tilde\psi_{k\ep}$ is uniformly Lipschitz continuous
and uniformly bounded, hence Ascoli-Arzel\`a's theorem gives the existence
of a Lipschitz function $\tilde\psi_k:D\lra\R$ that is the uniform limit
of $\tilde\psi_{k\ep}$ on compact sets of $B_1$ and clearly extends to the closure $D$.
Furthermore, $d\tilde\psi_k=f_k$ as distributions, hence taking into account that
$\tilde\psi_k$ is Lipschitz, it follows that $d\tilde\psi_k=f_k$ a.e. in $D$.
Thus, defining $\psi_k:D\lra\R$ as $\psi_k=\tilde\psi_k+z_k(1,0)-\tilde\psi_k(1,0)$
for every $k=1,\ldots,s$, we get
\[
\psi_k(1,0)=z_k(1,0)\quad\mbox{and}\quad d\psi_k=\ph^*\theta^k\quad\mbox{a.e. in $D$}.
\]
As a result, the function $T:D\lra V_2$ given by $T=\sum_{k=1}^s\psi_k\,Z_k$
satisfies the condition $T(1,0)=\sum_{k=1}^sz_k(1,0)\,Z_k$ along with
the contact equations \eqref{contZ}.
The last step is to show the validity of the boundary condition
$T_{|S^1}=\sum_{k=1}^sz_k\,Z_k$.

To do this, we recall that the mapping
\[
\Phi:D\lra\A\cl_{\cn}^n,\quad \Phi=\exp\circ(A+T)
\]
is Lipschitz continuous. In particular, $\Phi_{|S^1}:S^1\lra\A\cl_{\cn}^n$ and 
$\Gamma:S^1\lra\A\cl_{\cn}^n$ are Lipschitz curves with $\Phi(1,0)=\Gamma(1,0)=e$\,.
Then these curves must coincide by uniqueness of the horizontal lifting.
In fact, we have
\[
\left\{\begin{array}{l}
\Phi(e^{2\pi i t})=\exp\big(a(t)+\tilde b(t)\big) \\
\Gamma\big(e^{2\pi it}\big)=\exp\big(a(t)+b(t)\big)
\end{array}\right.\,,
\]
where $b(0)=\tilde b(0)=0$. By Theorem~\ref{CharLip}, since 
both $\Phi_{|S^1}$ and $\Gamma$ are Lipschitz, then both $b$ and $\tilde b$ 
satisfy the contact equations
\[
\dot{\tilde b}=\frac{1}{2}[a,\dot a]=\dot b\quad\mbox{a.e. in $[0,1]$}
\]
and clearly $b=\tilde b$. This shows that $\Phi_{|S^1}=\Gamma$ and leads us to the conclusion. $\Box$
\begin{Rem}{\rm
From definition of 1-connectedness, one easily observes
that Theorem~\ref{1Lipcon} exactly coincides with Theorem~\ref{Allcock1Lipcon}.
}\end{Rem}
\section{Quadratic isoperimetric inequalities}\label{Sectquadisop}
In this section, we prove the validity of quadratic isoperimetric inequalities
in Allcock groups. We will follow conventions and notation of  
Section~\ref{LipExtSect}.

Let $g$ be a left invariant Riemannian metric defined on $\A\cl^n$
such that the fixed basis $(X_{ij},Z_k)_{i,j,k}$ of $\cA\cl^n$ is orthonormal.
Denote by $\rho$ the associated Carnot-Carath\'eodory distance defined
on $\A\cl^n$ as 
\[
\rho_0(x,y)=\inf_{\Gamma:x\lra y,\,T\Gamma\subset H\A\cl^n}
\int_0^1 \sqrt{g\big(\Gamma(t)\big)\big(\Gamma'(t),\Gamma'(t)\big)}\,dt
\]
where $\Gamma:[0,1]\lra\A\cl^n$. Notice that if
$\Gamma(t)=\exp\big(\gamma_1(t)+\gamma_2(t)\big)$,
$\gamma_i(t)\in V_i$ and also
$\Gamma'(t)=\sum_{ij}a_{ij}(t)\,X_{ij}\in L^1\left((0,1),V_1\right)$,
then
\[
\int_0^1 \sqrt{g\big(\Gamma(t)\big)\big(\Gamma'(t),\Gamma'(t)\big)}\,dt=
\int_0^1 |a(t)|\,dt=\int_0^1|\gamma_1'(t)|\,dt.
\]
The last equality follows from \eqref{dexp}, taking into account
that $\Gamma$ is horizontal. The symbol $|\cdot|$ above 
also denotes the Hilbert norm in $\cA\cl^n$ that makes $(X_{ij},Z_k)$ orthonormal.
Abusing notation, we will use the same symbol to define a norm on $\A\cl^n$ as follows
\[
|x-y|:=|\exp^{-1}(x)-\exp^{-1}(y)|\,,
\]
for every $x,y\in\A\cl^n$.
We denote by $\varrho$ the Riemannian distance associated to $g$
and notice that $\varrho\leq\rho_0$.
\begin{Pro}\label{Vcomm}
Let $\V$ be a horizontal subgroup of $\A\cl^n$.
Then for every $x,y\in\V$, we have $\rho_0(x,y)=\varrho(x,y)=|x-y|$.
\end{Pro}
{\sc Proof.}
Let $x=\exp\xi$ and $y=\exp\eta$ and consider the curve $\widetilde\Gamma(t)
=\exp\big(\xi+t(\eta-\xi)\big)$,
where $\xi,\eta\in V$, where $V$ is the Lie algebra of $\V\subset\exp(V_1)$.
Then $[\xi,\eta]=0$ and
$\widetilde\Gamma$ is horizontal. Furthermore,
\[
\int_0^1 \sqrt{g\big(\tilde\Gamma(t)\big)\big(\tilde\Gamma'(t),
\tilde\Gamma'(t)\big)}\,dt
=\int_0^1|\widetilde\gamma_1'(t)|\,dt=|\xi-\eta|=|x-y|.
\]
If $\Gamma=\exp\circ(\gamma_1+\gamma_2)$
is any absolutely continuous curve connecting $x$ and $y$, then
$\gamma_1(0)=\xi$ and $\gamma_1(1)=\eta$ and clearly
\[
\int_0^1 \sqrt{g\big(\Gamma(t)\big)\big(\Gamma'(t),\Gamma'(t)\big)}\,dt
=\int_0^1|\gamma_1'(t)|\,dt\geq|\xi-\eta|=|x-y|.
\]
This concludes the proof. $\Box$
\begin{Pro}\label{lengthcurve}
Let $\Gamma:S^1\lra\A\cl^n$ be a Lipschitz mapping defined as
\[
\Gamma=\exp\circ\Big(\sum_{ij}c_{ij}\,X_{ij}+\sum_{k}z_k\,Z_k\Big)\,.
\]
Then we have $\mbox{\rm length}_\varrho\Gamma
=\mbox{\rm length}_{\rho_0}\Gamma=\mbox{\rm length}_{|\cdot|}(c)$.
\end{Pro}
{\sc Proof.}
Using the sub-Riemannian area formula, one gets
\[
\int_a^b\varrho\big(\exp(\dot\gamma_1(t))\big)\,dt=\int_a^b\rho_0\big(\exp(\dot\gamma_1(t))\big)\,dt
=\int_\G N(\Gamma,y)\,d\mathcal H^1(y)\,,
\]
then Theorem~2.10.13 of \cite{Fed} gives
\[
\int_a^b\varrho\big(\exp(\dot\gamma_1(t))\big)\,dt=\int_a^b\rho_0\big(\exp(\dot\gamma_1(t))\big)\,dt
=V_a^b\Gamma=\mbox{\rm length}_{\rho_0}(\Gamma)
=\mbox{\rm length}_{\varrho}(\Gamma)\,.
\]
Taking into account Proposition~\ref{Vcomm}, it follows that
$\rho_0\big(\exp(\dot\gamma_1(t))\big)=|\dot\gamma_1(t)|=|\dot c(t)|$. $\Box$
\begin{Def}{\rm
Let $L:\R^2\lra\A\cl^n$ be an h-homomorphism. Then the {\em jacobian}
of $L$ is given by
\[
J(L)=\frac{\mathcal H_{\rho_0}^2\big(L(A)\big)}{\mathcal H_{|\cdot|}^2(A)}
\]
where $A$ is any set of positive measure in the plane.
}\end{Def}
The previous definition of jacobian has been introduced in \cite{Mag}
for h-homomorphisms of stratified groups. Notice that it does not depend
on the choice of the set $A$.
\begin{Rem}\label{LhAl}{\rm 
Let $L:\R^2\lra\A\cl^n$ be an h-homomorphism. Then $L(\R^2)=\V$
is a horizontal subgroup of $\A\cl^n$ and Proposition~\ref{Vcomm} implies that
\[
\mathcal H_{\rho_0}^2\big(L(A)\big)=\mathcal H_{|\cdot|}^2\big(L(A)\big)\,,
\]
then the classical area formula yields
\[
J(L)=|L_1\wedge L_2|
\]
where $|L_1\wedge L_2|$ is the classical Euclidean jacobian of $L$
and $L_i=L(e_i)$.
}\end{Rem}
\begin{Pro}\label{srareapro}
Let $f:D\lra\A\cl^n$ be a Lipschitz mapping defined as
\[
f=\exp\circ\Big(\sum_{ij}\ph_{ij}\,X_{ij}+\sum_{k}\psi_k\,Z_k\Big)\,.
\]
Let $N(f,y)=\mathcal H^0\big(f^{-1}(y)\big)$ be the 
multiplicity function. Then we have
\begin{eqnarray}\label{srarea}
\int_{\A\cl^n}N(f,y)\,\mathcal H^2_{\rho_0}(y)
=\int_D\,J\big(Df(x)\big)\,dx=\int_D|\der_{x_1}\ph\wedge\der_{x_2}\ph|\,dx\,.
\end{eqnarray}
\end{Pro}
{\sc Proof.}
Since the Pansu differential $Df(x):\R^2\lra\A\cl^n$ is an h-homomorphism,
it has the following matrix representation
\[
Df(x)=\left(\begin{array}{c} 
 \nabla\ph_{11}  \\
 \nabla\ph_{12} \\
 \nabla\ph_{13} \\
  \vdots            \\
 \nabla\ph_{nm} \\
    0 \\
  \vdots \\
    0
\end{array}\right)\,.
\]
In view of Remark~\ref{LhAl}, the sub-Riemannian area formula of
\cite{Mag} concludes the proof. $\Box$
\begin{The}\label{projquadisoperthe}
Let $\A\cl^n$ be an Allcock group and let $\omega$ and $\theta$ be
the associated forms on $\R^{mn}$ with respect to fixed graded
coordinates, where $n\geq2.$ 
Then there exists a geometric constant $K>0$ such that
for every Lipschitz loop $c:S^1\lra\R^{mn}$, with $\int_c\theta=0$ and
$c(1,0)=0$ one can find a Lipschitz extension $\ph:D\lra\R^{mn}$
such that $\ph^*\omega=0$ a.e. in $D$ and
\begin{equation}\label{projquadisoper}
\mathcal H_{|\cdot|}^2\big(\ph(D)\big)\leq
\int_D|\der_{x_1}\ph\wedge\der_{x_2}\ph|\,dx\leq K\;
\mbox{\rm length}_{|\cdot|}(c)^2\,,
\end{equation}
where $\mbox{\rm length}_{|\cdot|}(c)$ is the length of $c$
with respect to the Euclidean norm $|\cdot|$ in $\R^{mn}$.
\end{The}
{\sc Proof.}
We will continue the argument used in the proof of Theorem~\ref{allcock},
exploiting the same notation. We first recall that
\[
\mbox{\rm length}_{|\cdot|}(c)=\int_0^1|\alpha'
(t)|\,dt
\]
is the length of $c:S^1\lra\R^{mn}$, where $\alpha(t)=c(e^{i2\pi t})$.
From the proof of Theorem~\ref{allcock}, we recall the definition of 
$\ph:D\lra\R^{mn}$ as
\[
\ph\big(\rho\,e^{i2\pi\theta}\big)=\left\{\begin{array}{ll}
H\big(\theta,2(1-\rho)\big)  &
 \mbox{if}\; 1/2\leq \rho\leq 1 \\
q_0  & \mbox{if}\; 0\leq \rho\leq1/2 \end{array}\right.\,,
\]
where $H:[0,1]^2\lra\R^{mn}$ is given by
\[
H(\tau,t)=\Gamma_{k+1}(\tau,5t-k)  
\quad\mbox{if}\quad \frac{k}{5}\leq t<\frac{k+1}{5}
\quad\mbox{and}\quad k=0,1,2,3,4\,.
\]
If we set $\tilde\ph(\rho,\theta)=\ph(\rho\,e^{i2\pi\theta})$, then we have
\[
\int_D\,|\ph_{x_1}\wedge\ph_{x_2}|\,dx=\int_{1/2}^1\int_0^1J\tilde\ph(\rho,\theta)
\,d\theta\,d\rho\,,
\]
where $|v\wedge w|$ is the Hilbert norm on 2-vectors of $\Lambda_2(\R^{mn})$
with respect to the canonical basis.
This is a consequence of the change of variable
$\phi(\rho,\theta)=\rho\,e^{i2\pi \theta}$ and the fact that
$|\phi_\rho\wedge\phi_\theta|=\rho$. It follows that
\[
\int_DJ\ph(x)\,dx=2\int_{1/2}^1\int_0^1
|H_\theta\big(\theta,2(1-\rho)\big)
\wedge H_\rho\big(\theta,2(1-\rho)\big)|\,d\theta\,d\rho\,.
\]
Taking into account that
\[
H\big(\theta,2(1-\rho)\big)=\Gamma_{k+1}\big(\theta,10(1-\rho)-k\big)  
\quad\mbox{if}\quad 1-\frac{k+1}{10}<\rho\leq1-\frac{k}{10}
\quad\mbox{and}\quad k=0,1,2,3,4\,,
\]
a simple change of variable yields
\begin{eqnarray}\label{areaDJph}
\int_D|\ph_{x_1}\wedge\ph_{x_2}|\,dx
=\sum_{k=0}^4\int_0^1\int_0^1|(\Gamma_{k+1})_\tau\big(\theta,t\big)
\wedge (\Gamma_{k+1})_t\big(\theta,t)|\,dt\,d\theta\,.
\end{eqnarray}
Now, we use the explicit formulas of $\Gamma_j$ given in the proof of
Theorem~\ref{allcock}.
From definition of $\Gamma_1$ it is obvious that $|(\Gamma_1)_\tau\wedge(\Gamma_1)_t|=0$ a.e. 
in the unit square $[0,1]^2$, denoted by $Q$. 
Simple computations yield the following estimates
\begin{eqnarray*}
\int_Q|(\Gamma_2)_\tau\wedge(\Gamma_2)_t|d\tau dt&=&
2\int_Q|\dot\beta_1(2\tau-1+t)\wedge\dot\beta_2(2\tau-1)|d\tau dt\\
&\leq&2\int_0^1\left(\int_0^1|\dot\beta_1(2\tau-1+t)|dt\right)|\dot\beta_2(2\tau-1)|d\tau\\
&\leq&2\int_\R\left(\int_\R|\tilde\alpha'(t)|dt\right)
|\tilde\alpha'(2\tau-1)|d\tau=2\left(\int_0^1|\alpha'(t)|dt\right)^2\,.
\end{eqnarray*}
The area contributed by $\Gamma_3$ is given by
\begin{eqnarray*}
\int_Q|(\Gamma_3)_\tau\wedge(\Gamma_3)_t|d\tau dt&=&\pi\int_Q\left|\left(\cos\big(\pi
t/2\big)\dot\beta_1(2\tau)+\sin\big(\pi t/2\big)
\dot{\overline{\beta_1}}(2\tau)+\dot\beta_2(2\tau-1)\right)\right.\\
&&\wedge\left(-\sin\big(\pi t/2\big)\beta_1(2\tau) +\cos\big(\pi
t/2\big)\overline{\beta_1}(2\tau)\right)\Big|d\tau dt\\
&\leq&\pi\int_Q\left|\left(\cos\big(\pi t/2\big)\dot\beta_1(2\tau)+\sin\big(\pi t/2\big)
\dot{\overline{\beta_1}}(2\tau)+\dot\beta_2(2\tau-1)\right)\right|\\
&&\left|\left(-\sin\big(\pi t/2\big)\beta_1(2\tau) +\cos\big(\pi
t/2\big)\overline{\beta_1}(2\tau)\right)\right|d\tau dt\\
&\leq&2\pi\int_0^1\big(2|\dot\beta_1(2\tau)|
+|\dot\beta_2(2\tau-1)|\big)\,|\beta_1(2\tau)|\,d\tau\\
&\leq&\pi\int_0^1|\dot\alpha(t)|dt\;
\int_0^2\big(2|\tilde\alpha'(t)|+|\tilde\alpha'(t-1)|\big)\,dt\\
&=&3\pi\left(\int_0^1|\dot\alpha(t)|dt\right)^2\,.
\end{eqnarray*}
Concerning the isotropic homotopy $\Gamma_4$, we get
\begin{eqnarray*}
\int_Q|(\Gamma_4)_\tau\wedge(\Gamma_4)_t|d\tau dt&=
&\frac{\pi}{2}\int_Q\Big|\left(\cos\big(\pi t/2\big)(\overline\beta_1\beta_2)'(\tau)
+\sin\big(\pi t/2\big)q'(\tau)\right)\\
&&\wedge\left(-\sin\big(\pi t/2\big)(\beta_1\beta_2)(\tau)+
\cos\big(\pi t/2\big)q(\tau)\right)\Big|d\tau dt\\
&\leq&\frac{\pi}{2}\int_0^1\big(|(\overline\beta_1\beta_2)'(\tau)|+
|q'(\tau)|\big)\,
\big(|(\overline\beta_1\beta_2)(\tau)|+|q(\tau)|\big)\,d\tau
\end{eqnarray*}
By definition of $\overline\beta_1\beta_2$, one easily checks that
\begin{equation}\label{maxbeta12}
\max_{[0,1]}|(\overline\beta_1\beta_2)|\leq2\int_0^1|\dot\alpha(t)|\,dt
\quad\mbox{and}\quad
|(\overline\beta_1\beta_2)'(t)|\leq2\big(|\tilde\alpha'(2t)|
+|\tilde\alpha'(2t-1)|\big)\,.
\end{equation}
The curve $q(t)$ is given by the proof of Theorem~\ref{allcock},
where it has been obtained applying Definition~\ref{suriso} with
$\lambda=\mbox{\rm length}_{|\cdot|}(c)$. It follows that there exists $C>0$
such that $|q(0)|\leq C\,\lambda$. Thus, we get
\[
|q(t)|\leq C\lambda+\frac{C}{\lambda}\int_0^t|\sigma(t)|\,dt\,.
\]
Recall that
\[
\sigma(t)=\sum_{k=1}^s\omega^k\left((\overline\beta_1\beta_2)(t)\wedge
(\overline\beta_1\beta_2)'(t)\right)\,Z_k\,,
\]
then bilinearity of $\omega^k$'s yields $C_1>0$ such that
\[
|\sigma(t)|\leq C_1
\,|(\overline\beta_1\beta_2)(t)|\,|(\overline\beta_1\beta_2)'(t)|,
\]
hence in view of \eqref{maxbeta12}, it follows that
\begin{equation}\label{estsigma}
|\sigma(t)|\leq 4C_1\,
\lambda\,\big(|\tilde\alpha'(2t)|+|\tilde\alpha'(2t-1)|\big)\,.
\end{equation}
As a result, we get
\[
\max_{t\in[0,1]}|q(t)|\leq C\,\lambda\,(1+4C_1)\,,
\]
that implies
\begin{eqnarray*}
\int_Q|(\Gamma_4)_\tau\wedge(\Gamma_4)_t|d\tau dt&\leq&
\lambda\;\frac{\pi}{2}\;\big(2+C+4C_1C\big)
\int_0^1\big(|(\overline\beta_1\beta_2)'(\tau)|+
|q'(\tau)|\big)\;d\tau\,.
\end{eqnarray*}
Taking into account that
\[
\int_0^1|\dot q(t)|\,dt\leq\frac{C}{\lambda}\int_0^1|\sigma(t)|\,dt
\leq 2\,C\,C_1\int_0^1
|(\overline{\beta}_1\beta_2)'(t)|\leq4\,C\,C_1\,\lambda\,,
\]
we are then lead to the following
\begin{eqnarray*}
\int_Q|(\Gamma_4)_\tau\wedge(\Gamma_4)_t|d\tau dt&\leq&
\pi\;\big(2+C+4C_1C\big)\,\big(1+2CC_1\big)\,\lambda^2\,.
\end{eqnarray*}
Finally, by definition of Allcock group and applying Definition~\ref{suriso}, we have
\begin{eqnarray*}
\int_Q|(\Gamma_5)_\tau\wedge(\Gamma_5)_t|d\tau dt&\leq&\,C\,
\left(\int_0^1|\dot q(t)|\,dt\right)^2\leq 16\,C^3\,C_1^2\;\lambda^2\,.
\end{eqnarray*}
Joining the previous estimates with \eqref{areaDJph},
estimate \eqref{projquadisoper} follows. $\Box$
\begin{Rem}{\rm
The previous theorem could be seen as a
completion of Theorem~\ref{allcock}, where we have fixed our attention
on the area enclosed by the extension $\ph$.
}\end{Rem}
{\sc Proof of Theorem~\ref{quadisop}.}
Up to left translation, we can assume that $\Gamma(1,0)$ coincides
with the unit element. Thus, we define 
\[
\Gamma=\exp\circ\left(\sum_{i=1}^n\sum_{j=1}^nc_{ij}\,X_{ij}+
\sum_{k=1}^sz_k\,Z_k\right)\,,
\]
where $c:S^1\lra\R^{mn}$. Arguing as in the beginning of the proof of  Theorem~\ref{1Lipcon}, it follows that $\int_c\theta=0$
and $c(1,0)=0$. Then we apply Theorem~\ref{projquadisoperthe},
getting a Lipschitz extension $\ph:D\lra\R^{mn}$ such that
\begin{equation}
\int_D|\der_{x_1}\ph\wedge\der_{x_2}\ph|\,dx\leq K\;
\mbox{\rm length}_{|\cdot|}(c)^2\,,
\end{equation}
Finally, Proposition~\ref{lengthcurve} and Proposition~\ref{srareapro}
lead us to the conclusion. $\Box$
\begin{Exa}\label{exaheis}{\rm
Let us assume the setting in the proof of Theorem~\ref{quadisop},
where the Allcock group is the 7-dimensional Heisenberg group
$\H^2$. According to this proposition, we set $\phi=(\ph,\psi)$.
We equip $\H^2$ with suitable graded coordinates $(x_1,\ldots,x_5)$
such that
\begin{eqnarray}\label{nablaheis2}
\nabla\phi(x)=\left(\begin{array}{c} 
 \nabla\ph_1  \\
 \nabla\ph_2 \\
 \nabla\ph_3 \\
 \nabla\ph_4  \\
 \ph_1\nabla\ph_3-\ph_3\nabla\ph_1+\ph_2\nabla\ph_3-\ph_3\nabla\ph_2 \\
\end{array}\right)\,,
\end{eqnarray}
where the contact equations imply that $\nabla\psi$ is equal
to the last row of \eqref{nablaheis2}.
Then a simple computation yields
\[
|\der_{x_1}\phi\wedge\der_{x_2}\phi|\leq
\sqrt{1+3|\ph|^2}\;|\der_{x_1}\ph\wedge\der_{x_2}\ph|\,.
\]
Now, if $\ph$ is the extension provided by Theorem~\ref{projquadisoperthe}
and $f=\exp\left(\sum_{j=1}^4\ph_j\,X_j+\psi\,X_5\right)$
is the corresponding Lipschitz mapping, then
\[
\mathcal H^2_{|\cdot|}\left(f(D)\right)\leq C
\left(1+2\max_{D}|\ph|\right)\,\mathcal \mbox{\rm length}_{|\cdot|}(c)^2
\]
for a suitable geometric constant $C>0$.}\end{Exa}

\end{document}